\documentclass[11pt, lettersize, notitlepage, leqno, footskip=0.6cm]{article} 
\usepackage{amssymb}
\usepackage{amsthm}
\usepackage{amsmath}	
\usepackage{graphicx}
\usepackage{amscd}

\usepackage{thmtools}
\usepackage[none]{hyphenat}	

\newcommand{\bc}{\mathbb C}
\newcommand{\bF}{\mathbb F}

\newcommand{\bz}{\mathbb Z}

\newcommand{\bq}{\mathbb Q}
\newcommand{\br}{\mathbb R}
\newcommand{\bS}{\mathbb S}

\newcommand{\bT}{\mathbf T}
\newcommand{\bTl}{\mathbf T_{{\bq_l}}}

\newcommand{\G}{\mathcal G}

\newcommand{\Gal}{\mathrm{Gal}}

\newcommand{\W}{\mathcal W}
\newcommand{\WA}{\mathcal{W}_{A_v}}

\newcommand{\zbar}{\overline {\mathbb{Z}}}
\newcommand{\qbar}{\overline {\mathbb{Q}}}

\newcommand{\Kbar}{\overline {K}}

\newcommand{\bark}{\overline {k}}

\newcommand{\bg}{\mathbb{G}}
\newcommand{\ev}{\mathbf{ev}}

\newcommand{\lcm}{\mathrm{lcm}}

\newcommand{\absq}{\mathrm{Gal}_{\bq}}
\newcommand{\absql}{\mathrm{Gal}_{\bq_l}}

\newcommand{\absf}{\mathrm{Gal}_F}
\newcommand{\absfv}{\mathrm{Gal}_{F_v}}

\newcommand{\absk}{\mathrm{Gal}_K}
\newcommand{\absl}{\mathrm{Gal}_L}

\newcommand{\Fr}{\mathrm{Fr}}

\newcommand{\Gm}{\mathbb{G}_m}

\newcommand{\la}{\langle}
\newcommand{\ra}{\rangle}

\newcommand{\lra}{\longrightarrow}
\newcommand{\xra}{\xrightarrow}
\newcommand{\hra}{\hookrightarrow}

\newcommand{\ti}{\tilde}
\newcommand{\wti}{\widetilde}
\newcommand{\mf}{\mathfrak}
\newcommand{\mc}{\mathcal}
\newcommand{\mb}{\mathbb}
 
\newcommand{\mr}{\mathrm}

\newcommand{\mfp}{\mathfrak{p}}
\newcommand{\tmfp}{\ti{\mf{p}}}

\newcommand{\setm}{\setminus}

\newcommand{\Br}{\mr{Br}}

\newcommand{\al}{\alpha}
\newcommand{\be}{\beta}
\newcommand{\ga}{\gamma}

\newcommand{\Gam}{\Gamma}

\newcommand{\Lamb}{\Lambda}

\newcommand{\si}{\sigma}

\newcommand{\Om}{\Omega}

\newcommand{\ov}{\overline}
\newcommand{\brho}{\ov{\rho}}

\newcommand{\ovp}{\ov{\pi}}

\newcommand{\sub}{\subseteq}

\newcommand{\bqW}{\bq(\W_{A_v})}
\newcommand{\bqWr}{\bq(\W_{A_v})^{^+}}
\newcommand{\GalW}{\Gal(\bq(\W_{A_v})/\bq)}
\newcommand{\GalWr}{\Gal(\bq(\W_{A_v})^+/\bq)}

\newcommand{\GA}{\mathbf{G}_A}
\newcommand{\VA}{\mathbf{V}_A}
\newcommand{\GAQ}{\mathbf{G}_{A,\qbar}}

\DeclareMathOperator{\MT}{MT}

\DeclareMathOperator{\Aut}{Aut}
\DeclareMathOperator{\Cent}{Cent}

\DeclareMathOperator{\End}{End}
\DeclareMathOperator{\GL}{GL}
\DeclareMathOperator{\GSU}{GSU}

\DeclareMathOperator{\Cor}{Cor}

\DeclareMathOperator{\Mat}{Mat}
\DeclareMathOperator{\PSL}{PSL}
\DeclareMathOperator{\SL}{SL}
\DeclareMathOperator{\SU}{SU}
\DeclareMathOperator{\GSp}{GSp}
\DeclareMathOperator{\Lie}{Lie}
\DeclareMathOperator{\rank}{rank\;}
\DeclareMathOperator{\id}{id}
\DeclareMathOperator{\inv}{inv}

\DeclareMathOperator{\ad}{ad}
\DeclareMathOperator{\der}{der}

\DeclareMathOperator{\charac}{char}
\DeclareMathOperator{\Res}{Res}

\DeclareMathOperator{\conn}{conn}

\newcommand{\Conj}{\mathbf{Conj}}
\newcommand{\Spec}{\mr{Spec}}

\newcommand{\Cl}{\mr{Cl}}

\newcommand{\ClGA}{\mr{Cl}_{\GA}'}

\newcommand{\ovPV}{\ov{P}_{_{\VA}}}

\newtheorem{Thm}{Theorem}[section]
\newtheorem{Prop}[Thm]{Proposition}
\newtheorem{Lem}[Thm]{Lemma}
\newtheorem{Corr}[Thm]{Corollary}

\newtheorem{Def}{Definition}[section]

\declaretheoremstyle[%
  spaceabove=-2pt,%
  spacebelow=8pt,%
  headfont=\normalfont\itshape,%
  postheadspace=1em,%
  qed=\qedsymbol%
]{mystyle} 
\declaretheorem[name={Proof},style=mystyle,unnumbered,
]{prf}

\numberwithin{equation}{section}

\title{Reductions of abelian varieties of generalized Mumford type}
\author{Steve Thakur}  
\date{\vspace{-3ex}}

\begin{document} 
\maketitle

\begin{abstract} We study the special fibers of a certain class of absolutely simple abelian varieties over number fields with endomorphism rings $\bz$ and possessing $l$-adic monodromy groups of the least possible rank. We also study the Dirichlet density of the places at which the possible reductions occur and confirm a special case of a broader conjecture for the splitting of reductions of abelian varieties over number fields. 
\end{abstract}

\begin{center}
\section{\fontsize{11}{11}\selectfont Introduction}
\end{center}
\vspace{-0.3cm}

For an absolutely simple abelian variety $A$ over a number field $F$, it is natural to study the reduction $A_v$ at each place $v$ of good reduction. In particular, it is a question of significant interest to determine the splitting of $A_v$ into simple abelian varieties up to isogeny. The following conjecture suggests that the endomorphism algebra $\End^0_{\qbar}(A)$ determines the splitting of $A_v$ at almost all places: 

\noindent \textbf{Conjecture 1.} ([MP08], [Zyw14]) \textit{Let $X$ be an absolutely simple abelian variety over a number field $F$. After replacing $F$ by a finite extension if necessary, there exists a density one set of places such that the reduction $X_v$ is isogenous to the $d$-th power of a simple abelian variety where $d^2$ is the dimension of the division algebra $\End^0_{\qbar}(X)$ over its center.}

Substantial progress towards this conjecture was made in the main theorem of [Zyw14] which settles this conjecture for abelian varieties for which the Mumford-Tate conjecture is true. For instance, this implies conjecture 1 for abelian varieties with complex multiplication. Furthermore, for an abelian variety $X$ of dimension $g$ over a number field with $\End_{\qbar}(X) = \bz$, ([Pin98], Theorem 5.13) states that unless $2g$ lies in the (rather thin) set\vspace{-0.2cm} $$\{n^{2k+1}:n \geq 2,\;k\geq 1 \}\bigcup \{{{4k+2}\choose {2k+1}} :k\geq 1\} ,$$ the inclusion $G_l \sub \GSp_{2g,\bq_l}$ is an equality and hence, the Mumford-Tate conjecture holds for $X$. The first integer in this exceptional set is $4$ and the dimension $4$ case remains open to this day. 

The primary goal of this article is to study Conjecture $1$ (and related questions) for abelian varieties \textit{of Mumford type} -- a particular class of abelian varieties with endomorphism rings $\bz$ for which the Mumford-Tate conjecture remains open. They furnish us with examples of abelian varieties that have endomorphism rings $\bz$ but have $l$-adic monodromy groups $G_l$ much smaller than $\GSp_{2g,\bq_l}$; in fact the rank of $G_l$ is the least afforded by Orr's inequality ([Orr15]) \vspace{-0.15cm}$$\rank G_l(A)\geq \log_2 d+2.$$ Since the Mumford-Tate conjecture remains open for this class of abelian varieties, it is important to make a distinction between what we call abelian varieties \textit{weakly of Mumford type} and those \textit{strongly of Mumford type} (see definitions 1.2 and 1.3). The former definition imposes a condition on the groups $G_l$ while the latter imposes the analogous condition on the Mumford-Tate group $\MT(A)$. (A justification for this terminology can be found in Proposition 3.7 where we verify that an abelian variety strongly of Mumford type is weakly of Mumford type.) 

By definition, for an abelian variety \textit{weakly} of Mumford type, the group $G_l$ has a one-dimensional center and its derived group $G_{l,\qbar_l}^{\der}$ is a product of copies of $\SL_{2,\qbar_l}$. Relying heavily on a combination of techniques from [Pin98], [No09] and [Zyw14], we deduce the following result which settles Conjecture 1 for this class of abelian varieties:
\vspace{0.1cm}
\begin{Thm} Let $A$ be an abelian variety over a number field $F$ such that:\\
$1.$ The $l$-adic monodromy groups $G_l$ are all connected.\\
$2.$ For all primes $l$ in a set of positive density, $\Lie(G_l)_{\qbar_l}\cong \mb{G}_{a,\qbar_l}\oplus \left( \bigoplus_{i=1}^{N}\mf{sl}_{2,\qbar_l} \right)$ with $\mf{sl}_{2,\qbar_l}^{N}$ acting on $V_l\otimes_{\bq_l}\qbar_l$ by the $N$-th external tensor power of the standard representation of $\mf{sl}_{2,\qbar_l}$. 

\noindent Then the reduction $A_v$ is absolutely simple away from a set of places of Dirichlet density zero.\end{Thm}

The first condition is not particularly restrictive since by [LP95], we can always enlarge the number field to ensure that the $G_l$ are connected. Furthermore, for any abelian variety with endomorphism ring $\bz$, $G_l$ has a one-dimensional center and $\Lie(G_l^{\der})_{\qbar_l}$ is  isomorphic to a power of a simple Lie algebra of type $A$, $B$, $C$ or $D$. So the second condition in Theorem 1.1 amounts to $\End(A_{\qbar})$ being $\bz$ and this simple Lie algebra being of type $A_1$ for a set of rational primes of positive density. Table 4.2 of [Pin98] then implies that $V_l\otimes_{\bq_l}\qbar_l$ must necessarily be the $N$-th external tensor product of the standard representation of $\mf{sl}_{2,\qbar_l}$.

We call such abelian varieties \textit{weakly} of Mumford type so as to distinguish them from abelian varieties whose Mumford-Tate groups fulfill the analogous property. 

In the dimension $16$ (equivalently, rank $6$) case, we list the possible isogeny types for the reduction $A_v$ (Proposition 4.8) and briefly describe how the same techniques may be used for higher dimensions when the common rank of the derived subgroups of the $l$-adic monodromy groups is a prime. In the case of good ordinary reduction - which occurs with density one - we show that the reduction is either simple or has an ordinary elliptic curve as a simple component. We also show that both types of reduction occur.

\underline{A brief outline} In section 2, we list some of the fundamental results about Mumford-Tate groups and $l$-adic monodromy groups that we shall repeatedly need. In section 3, we prove a few short propositions that will be necessary in sections 4 and 5. In section 4, we describe the splitting of the reduction $A_v$ for an abelian variety weakly of Mumford type of rank $6$. We will follow much of the approach from [No00] but will need some more results about abelian varieties over finite fields and cyclic lattices (appendix 6.2). In section 5, we give a proof of the density result. The key result in the proof will be a variant of a compatibility theorem from [No09] that replaces the Mumford-Tate group by a certain reductive $\bq$-group constructed in [Pin98]. We use the appendix for fundamental results that are original but of an elementary nature.
\vspace{-0.3cm}
\subsection{\fontsize{11}{11}\selectfont  Notations and Terminology}
\vspace{-0.3cm}
For an abelian variety $A$ over a number field $F$ and a prime $l$, the \textit{$l$-adic monodromy subgroup} $G_l$ is the Zariski closure of the image of $\rho_{_l}:\absf\lra \GL(V_l(A))$ and $G_l^0$ is the connected component of the identity. The image $\rho_l(\absf)$ is a compact $l$-adic analytic subgroup of $G_l(\bq_l)$ which is of finite index in $G_l(\bq_l)$ ([Bog80]). We denote the Zariski closure of $\rho_{_l}(\absf)$ in $\GL(T_l(A))$ by $\mc{G}_l$. It is a group scheme over $\bz_l$ with generic fiber $G_l$.

Let $g$ be the dimension of $A$. The Betti cohomology group $H^1(A(\bc),\bq)$ is a $2g$-dimensional vector space over $\bq$ endowed with a decomposition $V\otimes_{\bq}\bc = V^{1,0}\oplus V^{0,1}$ such that $\overline{V^{1,0}} = V^{0,1}$. Let $\mu_{\infty}:\bg_{m,\bc}\lra \GL_{2g,\bc}$ be the cocharacter through which any $z\in \bc^{*}$ acts by multiplication by $z$ on $V^{1,0}$ and trivially on $V^{0,1}$. The \textit{Mumford-Tate group} $\MT(A)$ is the unique smallest $\bq$-algebraic subgroup of $\GL_{2g,\bq}$ such that $\mu_{\infty}$ factors through $\MT(A)\times_{\bq} \bc$.

The following long-standing conjecture suggests an intrinsic relation between the two notions:

\noindent\textbf{Conjecture} (Mumford-Tate) \textit{For an abelian variety $A$ over a number field and a prime $l$, $\MT(A)\times_{\bq}\bq_l = G_l^0$}.

While the conjecture remains open, several cases have been proven. In particular, if the endomorphism ring is trivial, the smallest dimension for which the conjecture remains open is $4$ (following the results of [Pin98]). Furthermore, the Mumford-Tate conjecture for abelian varieties over number fields implies Conjecture 1 ([Zyw14]) and hence, is closely related to the simple decomposition of the special fibers.

In this article, we study a class of abelian varieties with $l$-adic monodromy groups of low rank. The following is a straightforward generalization of the construction in [Mum69]. It is well-defined over an arbitrary field of characteristic zero. But throughout this article, we will only encounter cases where the field $K$ is either a number field or a local field.

\begin{Def}\normalfont Let $K$ be a field of characteristic zero, $G$ an algebraic group over $K$ and $V$ a faithful $K$-linear representation of $G$. Then $(G,V)$ is of \textit{Mumford type of rank $N+1$} if:\\
-$\Lie(G)$ has a one-dimensional center\\
-$\Lie(G)_{\Kbar}\cong \bg_{a,\Kbar}\oplus \underbrace{\mf{sl}_{2,\Kbar}\oplus\cdots \oplus \mf{sl}_{2,\Kbar}}_{_N}$\\
-$\mf{sl}_{2,\Kbar}^{N}$ acts on $V_{\Kbar}$ through the $N$-th (external) tensor power of the standard representation of $\mf{sl}_{2,\Kbar}$.\end{Def}
\vspace{-0.15cm}
If $G$ is further assumed to be connected, there is a central isogeny $\Phi:\ti{G}\lra G$ where $\ti{G}$ is an algebraic group over $K$ such that $\ti{G}\cong \bg_{m,K}\times \ti{G}^{\der}$ and $\ti{G}_{\Kbar}^{\der}\cong \SL_{2,\Kbar}^{N}$. We introduce the following notions so as to avoid ambiguity:

\begin{Def}\normalfont An abelian variety $A$ over a number field $F$ is \textit{weakly} of Mumford type of rank $N+1$ if $(G_l,V_l)$ is of Mumford type of rank $N+1$ for a set of rational primes $l$ of density one.\end{Def}
\vspace{-0.2cm}
We note that for such an abelian variety, $\mr{rank}(G_l(A)) = \log_2 (\dim A)+2$, meaning that Orr's inequality is an equality in this case. The rank of $G_l(A)$ is the same for all primes $l$ by a theorem of Serre's and hence, this equality also holds for the primes outside the density one set.

\begin{Def}\normalfont An abelian variety $A$ over a number field $F$ is \textit{strongly} of Mumford type of rank $N+1$ if $(\MT(A),H^1_B(X(\bc),\bq))$ is of Mumford type of rank $N+1$. \end{Def}
\vspace{-0.2cm}

We note that in both cases, the definitions imply that the abelian variety is of dimension $2^{N-1}$ and has endomorphism ring $\bz$. Furthermore, for an abelian variety weakly of Mumford type, $\mr{rank}(G_l(A)) = \log_2 (\dim A)+2$, meaning that Orr's inequality is an equality in this case.

Let $K$ be a totally real field of odd degree. Let $D$ a quaternion algebra over $K$ split at precisely one archimedean place with the corestriction $\Cor_{K/\bq}(D) =0$ in $\Br(\bq)$. Then an abelian variety corresponding to a closed point on the Shimura curve given by the datum $((D^*)^{\der},h)$ fulfills the conditions of Definition 1.3. We refer the reader to ([SZ12], Section 2) for details. We will see in appendix 6.1 that the converse is also true. 

An implication of ([Pin98], Theorem 5.15) is that an abelian variety \textit{strongly} of Mumford type is \textit{weakly} of Mumford type. The converse remains an open case of the Mumford-Tate conjecture (as far as we know). 

For an abelian variety $A$ with $\End(A_{\qbar}) = \bz$, we fix a triple $(\GA,\VA,\mb{S}_A)$ where:\\
- $\GA$ is a connected reductive group over $\bq$ with $\GA^{\der}$ $\bq$-simple\\
- $\VA$ is a faithful absolutely irreducible $\bq$-linear representation of $\GA$\\
- $\mb{S}_A$ is a density one set of rational primes such that\vspace{-0.1cm}$$(\GA,\VA)\otimes_{\bq}\bq_l\cong (G_l,V_l)\;\;\;\;\text{ for all } l\in \mathbb{S}_{A}.$$
The existence of this triple follows from [Pin98, Theorem 5.13]. Furthermore, if $A$ is weakly of Mumford type, then the pair $(\GA,\VA)$ is of Mumford type. If $A$ is \textit{strongly} of Mumford type, meaning that it arises as a closed point on a Mumford curve, the Mumford-Tate conjecture is true for $A$ by [Pin98]. So the group $\GA$ may be chosen to be the Mumford-Tate group $\MT(A)$ and $\VA$ the Betti cohomology group $H^1_B(A(\bc),\bq)$.
\vspace{-0.3cm}
\begin{center}
\section{\fontsize{11}{11}\selectfont Some background} 
\end{center}
\vspace{-0.3cm}
In this section, we state some well-known but fundamental results that we shall need repeatedly.

\begin{Thm}[Pin98] Let $X$ be an abelian variety over a number field $F$ such that $\End_{\qbar}(X) = \bz$ and the root system of each factor of $\MT(X)_{\qbar}$ has type $A_{2s-1}$ with $s\geq 1$ or $B_r$ with $r\geq 1$. Then the Mumford-Tate conjecture is true for $X$.\end{Thm}	
\vspace{-0.2cm}
In particular, if $X$ is strongly of Mumford type, the root system of every factor is $A_1$ and the Mumford-Tate conjecture is true for $X$. So $G_l(X) = \MT(X)\times_{\bq} \bq_l$ and 
\vspace{-0.1cm}$$\Lie(G_l(X))_{\qbar_l}\cong \Lie(\MT(X))\otimes_{\bq}\bq_l\cong \bg_{a,\qbar_l}\oplus \mathfrak{sl}_{2,\qbar_l}^{N}.$$

\begin{Thm} $\mr{[LP95]}\;$ For an abelian variety over a number field $F$, there exists a minimum finite extension $F_A^{\conn}$ where all the $G_l$ are connected. \end{Thm}
\vspace{-0.2cm}
This extension corresponds to the subgroup $\rho_{_l}^{-1}(G_l^0(\bq_l))\subseteq \absf$ which is independent of $l$. We will usually be replacing $F_A$ by $F_A^{\conn}$ throughout this article. This has the effect of replacing every $G_l$ by the connected component of the identity. The number field $F_A^{\conn}$ may be alternatively described as the intersection $$F_A^{\conn}:=\bigcap\limits_l F(A[l^{\infty}]). $$

\vspace{-0.2cm}

\begin{Thm} $\mr{(Faltings)}\;$ Let $X$ be an abelian variety over a number field $F$ such that $F_A^{\conn} = F$. Then we have the following:\\
\noindent $\mr{(i)}$ The centralizer of $G_l$ in $\End_{\bq_l}(V_l)$ is $\End(A)\otimes_{\bz} \bq_l$.\\
$\mr{(ii)}$ The group $G_l$ is reductive.\\
$\mr{(iii)}$ $\End^0(A_{F}) = \End^0(A_{\qbar})$\end{Thm}

\noindent The following result of Deligne gives one of the inclusions for the Mumford-Tate conjecture. 

\begin{Thm} $\mr{(Deligne)}\;$  For any prime $l$, $G_l^0\times_{\bq} \bq_l\subseteq \MT(A)\times_{\bq} \bq_l$.\end{Thm}

\vspace{-0.3cm}

\begin{Thm} $\mathrm{(Serre)}\;$ For any abelian variety $A$, the $l$-adic monodromy groups $G_l(A)$ are all of the same rank.\end{Thm}
\vspace{-0.2cm}
Henceforth, we shall refer to this common rank as \textit{the rank of the abelian variety}. Deligne's theorem on the first inclusion combined with Serre's rank independence theorem show that the Mumford-Tate conjecture reduces to the rank of the abelian variety being equal to the rank of the Mumford-Tate group. More precisely:

\begin{Thm} $\mathrm{(Serre)}\;$ If the inequality \vspace{-0.15cm}$$\rank G_p^{0}\leq  \rank \MT(A)$$ is an equality for any prime $p$, the Mumford-Tate conjecture is true for $A$.  \end{Thm}
\vspace{-0.2cm}
Furthermore, if we impose the additional condition that the endomorphism ring is trivial, as is the case with abelian varieties weakly of Mumford type, the following theorem from [Pin98] shows that there always exists a connected reductive group over $\bq$ that ``looks like" the Mumford-Tate group is expected to.

\begin{Thm} $\mr{([Pin98],\;Theorem\; 5.13)}$ Suppose $\End_{\qbar}(A) = \bz$. Then there exists a connected reductive group $\GA\sub \GL_{2g,\bq}$ over $\bq$ with a faithful absolutely irreducible $\bq$-linear representation $\VA$ such that $\GA^{\der}$ is $\bq$-simple and\vspace{-0.15cm}$$(\GA,\VA)\otimes \bq_l \cong (G_l,V_l)$$ for all rational primes $l$ in a set $\mb{S}_A$ of density one.\end{Thm}
\vspace{-0.2cm}
$\GA^{\der}$ being $\bq$-simple implies that $\Lie(\GA^{\der})$ is $\otimes$-isotypic. Furthermore, it is also known that the root systems of $\Lie(\GA^{\der})$ are of type $A$, $B$, $C$ or $D$ ([Pin98]).

\begin{Thm} $\mr{([LP95],\;Theorem\;3.2)}$ For an abelian variety $A$ over a number field $F$ such that $F= F_A^{\conn}$, $G_l$ is unramified $($quasi-split and split over some $\bq_{l^h}$$)$ for all but finitely many $l$. Furthermore, $\G_l$ is a reductive group scheme over $\bz_l$ for all but finitely many $l$.\end{Thm}
\vspace{-0.2cm}
(A reductive group over a local field is \textit{unramified} if it is quasi-split and split over some unramified extension.)
\vspace{0.05cm}

\noindent\underline{Strict Compatibility.} The $l$-adic representations attached to an abelian variety are \textit{strictly compatible} in the sense of Serre. A proof may be found in [Del74]. 

\begin{Thm} $\mr{[Del74]}$ For an abelian variety $A$, fix a finite set $\mathcal{S}$ of non-archimedean places such that $A$ has good reduction outside $\mathcal{S}$. Let $v\notin \mathcal{S}$, $l\neq \mathrm{char}(k_v)$\\
$1.$ $\rho_{_l}$ is unramified at $v$\\ 
$2.$ The characteristic polynomial of $\rho_{_l}(\mathrm{Fr}_v)$ has coefficients in $\bz$ and is independent of $l$.\end{Thm} 
\vspace{-0.2cm}
An immediate consequence is that the eigenvalues are in $\zbar$ and are independent of $l$. We denote this polynomial by $P_{A_v}(X)\in\bz[X]$, its zeros by $\WA$ and the multiplicative group they generate by $\Phi_{A_v}\sub \bq^*$. Note that a base change of the abelian variety to a field extension $\ti{F}/F$ and replacing $v$ by a place $\ti{v}$ lying over $v$ has the effect of replacing $\WA$ by $\WA^N:=\{\pi^N:\pi\in \WA \}$ where $N$ is the inertia degree of $v'$ over $v$.

Let $\bq(\W  _{A_v})^+$ be the compositum of all the totally real fields $\bq(\pi_{_i}+\ov{\pi}_{_i})$. Then $\Gal(\bq(\W  _{A_v})^+/\bq)\subseteq S_g$ and $\Gal(\bq(\W  _{A_v})/\bq(\W  _{A_v})^+)$ is abelian $2$-torsion of rank at most $g$. Thus, the Galois group of $P_{A_v}(X)$ is a subgroup of the semi-direct product $W_{2g}:=(\bz/2\bz)^{g}\rtimes S_{g}$, the ``generic Galois group" of the Galois closure of a CM field of degree $2g$. 

The Galois group is the full $W_{2g}$ away from a density zero set of places when the endomorphism ring is trivial and $G_l = \GSp_{\bq_l}$ [Cha97]. On the other hand, the Galois group is far smaller for abelian varieties weakly of Mumford type on account of the additive relations between the weights of the standard representation of $G_{l,\qbar_l}$ which induce multiplicative relations between the eigenvalues of the Frobenius.
\vspace{-0.4cm}
\begin{center}
\section{\fontsize{11}{11}\selectfont Some preliminary results}
\end{center}
\vspace{-0.2cm}
In this section, we discuss a few results we shall need in the subsequent sections. Some of them will be necessary to verify that the concepts from the preceding section are well-defined. 

\begin{Prop} Let $A$ be an abelian variety of dimension $2^{2n}$ over a number field $F$. Suppose there exists a prime $p$ in $\mathbb{S}_A$ such that $(G_p,V_p)$ is of Mumford type. Then $\End(A_{\qbar}) = \bz$\end{Prop}

\begin{prf} Passing to a finite extension if necessary, we may assume without loss of generality that all the $l$-adic monodromy groups are connected. Now, by Faltings, \vspace{-0.15cm}$$\End^0(A_{F'})\otimes_{\bq}\bq_p \cong \End_{\Gal_{F'}}(V_p) = \End_{G_p}(V_p) = \bq_p$$ for any finite extension $F'/F$. Thus, $\End^0(A_{\qbar}) = \bq$. \end{prf}
\vspace{-0.15cm}
We now verify that an abelian being weakly of Mumford type is equivalent to the hypothesis in Theorem 1.1.

\begin{Prop} Let $A$ be an abelian variety of dimension $2^{N-1}$ over a number field. Suppose the set of rational primes $p$ such $(G_p,V_p)$ is of Mumford type of rank $N+1$ has positive density. Then $(G_l,V_l)$ is of Mumford type for all $l\in \mathbb{S}_A$.\end{Prop}

\begin{prf} Since $\End(A_{\qbar}) = \bz$ by the preceding proposition, it follows from ([Pin98], Theorem 5.13) that there exists a connected reductive group $\GA$ over $\bq$ with a faithful $\bq$-linear representation $\VA$ such that $(\GA,\VA)\otimes_{\bq}\bq_l\cong (G_l,V_l)$ for all $l\in \mathbb{S}_A$. Since the set of rational primes described in the hypothesis has positive density, it intersects $\mb{S}_A$. So there exists a prime $p$ such that $(\GA,\VA)\otimes \bq_p\cong (G_p,V_p)$ and $(G_p,V_p)$ is of Mumford type. Hence, $(\GA,\VA)$ is of Mumford type of rank $N+1$ and the result follows.\end{prf}

\vspace{-0.2cm}

In the dimension $4$ case, $(G_p,V_p)$ being of Mumford type for a single prime $p$ implies the same for all primes. We refer the reader to [No00] for details. The proof follows from Serre's rank independence result and the classification of semisimple Lie algebras.

\begin{Lem} Let $G$ be a $\bq$-simple algebraic group. If $G$ is semisimple and the number of simple factors of $G_{\qbar}$ is a prime $p$, then the set of primes $l$ such that $G_{\bq_l}$ is $\bq_l$-simple is of positive density $\mr{(}$in particular, infinite$\mr{)}$.\end{Lem}
\begin{prf} The action of $\absq$ on the simple factors of $G_{\qbar}$ gives a map $\tau:\absq\rightarrow S_p$. Since $G$ is $\bq$-simple, the image of $\tau$ is a transitive subgroup of $S_p$ and hence, contains a $p$-cycle $c$. By the Chebotarev density theorem, the set of primes $l$ such that $c\in \tau(\absql)$ has positive density and in particular, is infinite. For any $l\in \mb{S}\cap \mc{P}$, $\absql$ acts transitively on the simple factors of $G_{\qbar_l}$ and hence, $G_{\bq_l}$ is $\bq_l$-simple. \end{prf}

\begin{Lem} Let $A$ be an abelian variety over number field such that $\End(A_{\qbar}) = \bz$. Suppose the number of simple factors of $G_{l,\qbar}^{\der}$ is a prime $p$ for all $l$ in a  set of density one. Then there exist infinitely many primes $l$ such that $G_{l}^{\der}$ is $\bq_l$-simple.\end{Lem}	
\begin{prf} From [Pin98], there exists a connected reductive group $\GA$ over $\bq$ with a faithful $\bq$-linear representation $\VA$ such that: 

\noindent$(\GA,\VA)\otimes \bq_l\cong (G_l,V_l)$ for all $l$ in a set $\bS_A$ of density one.\\
-$\GA^{\der}$ is $\bq$-simple.

Now $\Lie(\GA)$ has a one-dimensional center and $\Lie(\GA^{\der})$ is $\bq$-simple with $\Lie(\GAQ^{\der})$ having $p$ simple factors. The permutation action of $\absq$ on the simple factors of $\GAQ^{\der}$ gives a map $\tau:\absq\rightarrow S_{p}$. Since $\GA^{\der}$ is $\bq$-simple, it follows that the image is a transitive subgroup and since $p$ is a prime by assumption, the image contains a $p$-cycle $c$. But $c\in\tau(\absql)$ for a set $\mc{P}$ of rational primes of positive density by the Chebotarev density theorem. Since $\mb{S}_A$ is of density one, the intersection $\mb{S}_A\cap \mc{P}$ is also of positive density. For any prime $l\in \mb{S}_A\cap \mc{P}$, $G_{l}^{\der}$ is $\bq_l$-simple since the $p$-cycle acts transitively on the simple factors of $G_{l,\qbar}^{\der}$.\end{prf}
\vspace{-0.2cm}
We now use this lemma to show potential good reduction in the case when $\mr{rank}(G_l^{\der})$ is a prime. The key idea (following [No00]) is to use the fact that any abelian variety is de Rham everywhere, meaning that the inertia group acts unipotently on the Tate module of every different prime after passing to a finite extension, and subsequently show that the inertia group acts trivially on the Tate module.

\vspace{-0.1cm}
\begin{Prop} Suppose $N$ is a prime. Let $A$ be an abelian variety of rank $N+1$ weakly of Mumford type over a number field $F$. Then $A$ has potential good reduction everywhere.\end{Prop}

\begin{prf} Enlarging the field if necessary, we may assume without loss of generality that all the $G_l$ are connected and $A$ has semistable reduction everywhere. Since $N$ is assumed to be a prime, there exist infinitely many primes $l$ such that $G_l^{\der}$ is $\bq_l$-simple.

Let $v$ be a finite place of $F$ and let $I_{v}$ be the inertia group. Let $p$ be the rational prime lying under $v$. Choose a prime $l\neq p$ such that $G_l^{\der}$ is $\bq_l$-simple and let $\rho_{_l}$ be the $l$-adic representation associated to $A$. Since $A$ is semistable at $v$, the inertia group $I_{v}$ is rank $2$ unipotent on $V_l(A)$, meaning that $(\rho_{_l}(\sigma)-\id)^2 =0$ for every $\sigma\in I_{v}$ ([GRR72]). Suppose, by way of contradiction, that $\rho_{_l}(\sigma)\neq \id$ for some $\sigma\in I_v$. Now $\rho_{_l}(\sigma)\in G_l(\qbar)$ is the image by tensor product of an element $(y_{_1},\cdots,y_{_N})\in \SL_2(\qbar_l)$ with at least one $y_i\neq \id$. But since $G_l^{\der}$ is $\bq_l$-simple, it follows that $y_i\neq \id$ for any $i$. Conjugating if necessary, we may assume without loss of generality that all the $y_i$ are of the shape $$\begin{pmatrix}
1 & a_i\\ 
0 & 1\end{pmatrix}$$ with each $a_i\in \qbar_l\setm \{0\}$. So the index of nilpotency of $x-\id$ is larger than $2$, a contradiction. Hence, $I_{v}$ acts trivially on $V_l$ after passing to a finite extension. From the Neron-Ogg-Shafarevich criterion ([ST68]), it follows that $A$ has good reduction at $v$. \end{prf}
\vspace{-0.2cm}
For abelian varieties \textit{strongly} of Mumford type, [SZ12] have a proof of potential good reduction that does not assume $N$ being a prime. They do, however, assume that $v$ is prime to the discriminant of the totally real field of degree $N$ and to the discriminant of the quaternion algebra over that field whose group of norm one units yield the Shimura curve.

Furthermore, ([Pin98], Theorem 7.1) implies that the reduction is ordinary at most places.

\begin{Prop} For an abelian variety $A$ weakly of Mumford type over a field $F = F_A^{\conn}$, the places of good ordinary reduction have Dirichlet density $1$.\end{Prop}
\begin{prf} For any prime $l\in\mb{S}_A$, the root system of $G_{l,\qbar_l}^{\der}$ has type $A_1$ for a set of primes of density one. So the same is true for the root system of $\GAQ^{\der}$. In particular, the simple factors of $\GAQ^{\der}$ are not of type $C_r$ for any $r\geq 3$ and hence, $A$ has good ordinary reduction away from a density zero set of places by ([Pin98], Theorem 7.1).\end{prf}
\vspace{-0.2cm}
An alternate proof independent of the results of [Pin98] may be found in [No95]. The following proposition justifies the terminology in Definitions 1.2 and 1.3.

\begin{Prop} If an abelian variety $X$ is strongly of Mumford type, then it is weakly of Mumford type of the same rank.\end{Prop}
\begin{prf} We write $N = \rank \MT(X) -1 = \rank \MT(X)^{\der}$. By the hypothesis, the root system of every simple factor of $\MT(X)_{\qbar}$ has type $A_1$. So, by ([Pin98], Theorem 5.15), the Mumford-Tate conjecture is true for $X$. Hence, \vspace{-0.15cm}$$\Lie(G_{X,l})_{\qbar_l}\cong \Lie(\MT(X)\otimes_{\bq} \bq_l)_{\qbar_l} = \qbar_l\oplus \mf{sl}_{2,\qbar_l}^{N}.$$
Alternatively, note that if $X$ is strongly of Mumford type, then $\rank\MT(X) = N+1$. By [Orr15], we have the inequality $\mr{rank}(G_{X,l})\geq N+1$. Hence, $\mr{rank}(G_{X,l})\geq \rank \MT(X)$, which by Theorem 2.6, implies the Mumford-Tate conjecture for $X$. That the representation is the external tensor power of the standard representation of $\mf{sl}_{2,\qbar_l}$ follows from Table 4.2 of [Pin98].
\end{prf}

\begin{Prop} Let $G$ be a reductive group over a number field $K$. Let $T$ be a maximal torus of $G$. The set of primes $\mf{p}$ of $K$ such that $T_{K_{\mf{p}}}$ is split has positive density. \end{Prop}

\begin{prf} Let $L$ be a finite extension splitting $T$ and let $\wti{L}$ be the Galois closure over $K$. Let $\mc{S}(L|K)$ be the set of primes of $K$ split completely in $L$. For any $\mf{p}\in \mc{S}(L|K)$ and any prime $\tmfp$ of $L$ lying over $\mfp$ , $K_{\mf{p}} \cong L_{\tmfp}$ and hence, $T_{K_{\mf{p}}}$ is split. Since $\mc{S}(L|K)$ has density $\frac{1}{[\wti{L}:K]}$, the result follows.\end{prf}
\vspace{-0.1cm}
For the next two lemmas, let $A$ be an absolutely simple abelian variety over a number field $F$. Let $v$ be a place such that $\End(A_{\bark})$ is commutative.

\begin{Lem} Suppose $\End(A_v)$ is commutative. Then none of the quotients $\pi_{_1}\pi_{_2}^{-1}$ of the eigenvalues of $\Fr_v$ are roots of unity.\end{Lem}  
\begin{prf} Suppose $\Fr_v$ has two eigenvalues $\pi_{_1}$, $\pi_{_2}$ such that $\pi_{_1}\pi_{_2}^{-1} = \zeta_{_N}$ for some integer $N$. Choose a finite extension $F'/F$ and a place $\ti{v}$ lying over $v$ with inertia degree $N$ over $v$. Then $\W_{A_{\ti{v}}} = \W_{A_v}^N$ and hence, $P_{A_{\ti{v}}}(X)$ has the zero $\pi_{_1}^N$ with multiplicity at least two. In particular, $P_{A_{\ti{v}}}(X)$ is inseparable and hence, $\End(A_{\ti{v}})$ is non-commutative, a contradiction. \end{prf}

\begin{Lem} The reduction $A_v$ is absolutely simple if and only if $P_{A_v}(X)$ is irreducible.\end{Lem}
\begin{prf} Since $\End(A_{\bark})$ is commutative, the zeros of \vspace{-0.15cm}$$P_{A_v}^N(X):=\prod\limits_{\pi\in \WA}{X-\pi^N}$$ are distinct for any integer $N\geq 1$. So $A_v$ is absolutely simple if and only $\absq$ acts transitively on $\W_{A_v}^N = \W _{A_v\otimes k_{v^N}}$ for every integer $N$. Let $\si\in\absq$ and let $\pi_{_1},\pi_{_2}$ be any two eigenvalues of $\Fr_v$. Then $\si(\pi_{_1}^N) = \pi_{_2}^N \Longleftrightarrow \left(\frac{\si(\pi_{_1})}{\pi_{_2}}\right)^N =1 \Longleftrightarrow \\ \si(\pi_{_1}) = \pi_{_2}$ by the preceding lemma. So $\absq$ acts transitively on $\W_{A_v}^N$ if and only if it acts transitively on $\W_{A_v}$, which is equivalent to $P_{A_v}(X)$ being irreducible.\end{prf}

The lemmas 3.10 and 3.10 will help us detect places of absolutely simple reduction. In particular, as we will see in section 4, for an abelian variety weakly of Mumford type of rank $N+1$ with the integer $N$ a prime, absolutely simple reduction only occurs when $\End(A_{\bark})$ is commutative.
\vspace{-0.2cm}
\begin{center}
\section{\fontsize{11}{11}\selectfont Isogeny Classes}
\end{center}
\vspace{-0.2cm} 
We now study the isogeny classes of the special fibers of abelian varieties of dimension $16$ (equivalently, rank $6$) weakly of Mumford type. In particular, we determine the possible dimensions of the simple factors of the reductions and their endomorphism algebras. We will follow some of the techniques from [No00] but will need a few more techniques than in the dimension $4$ case. We also briefly describe how the techniques could be adapted to certain higher dimensional cases.
\vspace{0.05cm}

\subsection{\fontsize{11}{11}\selectfont The Galois group of the characteristic polynomial}
\vspace{-0.3cm}
\underline{Some Notation:} Let $A$ be an abelian variety of rank $N+1$ weakly of Mumford type over a number field $F$. Let $v$ be a place of good reduction and let $l$ be a prime other than the residue characteristic of $v$ such that $(G_l,V_l)$ is of Mumford type. Let $\rho_{_l}:\absf\lra \GL_{2g}(\bz_l)$ be the $l$-adic Galois representation and $\Fr_v$ a geometric Frobenius. Let $\pi_{_1},\cdots,\pi_{_{2g}}$ be the eigenvalues of $\rho_{_l}(\Fr_v)$. The eigenvalues lie in $\zbar$ and are independent of $l$. We will denote the set of eigenvalues by $\W_{A_v}$ and the multiplicative group they generate by $\Phi_{A_v}\sub\qbar^*$. The irreducible decomposition of the polynomial $P_{A_v}(X) = \prod\limits_{i=1}^{2g}(X-\pi_{_i})\in \bz[X]$ determines the simple decomposition of $A_v$. In particular, $P_{A_v}(X)$ is separable if and only if $\End^0(A_v)$ is a CM algebra. Furthermore, if $P_{A_v}(X)$ is separable, $A_v$ is simple if and only if $P_{A_v}(X)$ is irreducible. The field $\bq(\W_{A_v})$ is a Galois extension of $\bq$ and \vspace{-0.15cm}$$\Gal(\bq(\W_{A_v})/\bq)\sub\Aut\left(X^*(\bT)\right) = \{\pm 1\}^{N}\rtimes S_{_N}$$ where $\bT$ is a maximal torus of $\GA$. More precisely, \vspace{-0.15cm}$$\Gal(\bqW/\bqWr)\sub W(\GA,\bT) = \{\pm 1\}^{N}\;\;\;,\;\;\;\GalWr\sub \mr{Out}(\bT) = S_{_N}.$$
If $F'/F$ is a finite extension and $v'$ is a place of $F'$ lying over $v$, we have $\W_{A_{v'}} = \W_{A_{v}}^d$ where $d$ is the inertia degree of $v'$ over $v$. If the integer $N$ is a prime, $\GalW = \{\pm 1\}^j\rtimes H$ with $H$ a subgroup of $S_{_N}$ containing a $N$-cycle and $1\leq j\leq N$. Also, it is easy to show that $2^j\equiv 1\pmod {N}$ in this case (Proposition 4.10).

\vspace{-0.3cm}
\subsection{\fontsize{11}{11}\selectfont Weil Numbers and Isokummerian abelian varieties}
\vspace{-0.3cm}
We first state a few long-known results that will be useful for determining the possible isogeny classes of the simple components of the reductions. Proofs may be found in [Tat66] and [Wat69].

\underline{Notation:} For a Weil $q$-integer $\pi$ (up to conjugacy), let $B_{\pi}$ be the corresponding simple abelian variety over $\bF _q$ (up to isogeny). 

\begin{Lem} Let $\pi$ be a Weil $q$-integer and $B_{\pi}$ be the corresponding abelian variety over $\bF _{q}$. The following are equivalent:\\
$(1)$ $B_{\pi}$ is supersingular, meaning that the Newton slopes are all $\frac{1}{2}$\\
$(2)$ $\pi^N \in \bq$ for some integer $N$\\
$(3)$ $\End^0_{\overline{\bF}_p}(B_{\pi}) = \End^0_{\bF _{q^2}}(B_{\pi})\cong \Mat_{g}(\bq_{p,\infty})$ where $\bq_{p,\infty}$ is the quaternion algebra over $\bq$ ramified only at $p$ and $\infty$.\end{Lem}
\vspace{-0.1cm}
We may enlarge the field so that the size $q$ of the residue field is a square and none of the quotients $\frac{\pi}{\sqrt{q}}$ are non-trivial roots of unity. So the only possible Weil $q$-integer corresponding to a supersingular abelian variety would be $\sqrt{q}$. 
\vspace{-0.05cm}
\begin{Lem} For a simple abelian variety $B_{\pi}$ over $\bF _q$, the following are equivalent:\\
$1.$ $\End^0(B_{\pi})$ is commutative.\\
$2.$ The Weil number $\pi$ is of degree $2\dim B_{\pi}$.\end{Lem}

\vspace{-0.3cm}
\begin{Lem} For a simple abelian variety $B_{\pi}$ of dimension $g$ over $\bF _q$, the following are equivalent:\\
$1.$ $B_{\pi}$ is ordinary, meaning that the $p$-torsion group scheme is of rank $g$.\\
$2.$ $(\pi,\ov{\pi})=(1)$ in the ring of integers of the Galois closure of $\bq(\pi)$.\\
$2.'$ $(\pi+\ov{\pi},q)=(1)$ in the ring of integers of the Galois closure of $\bq(\pi)$.\\
$3.$ The Newton polygon is $g\times 0, g\times 1$\end{Lem}
\vspace{-0.2cm}
\underline{The group $\Phi_{A_v}$}. For an abelian variety $B$ over a finite field, we denote the $\absq$-stable multiplicative group generated by the eigenvalues of the Frobenius by $\Phi_B$. For a Weil number $\pi$, $\Phi_{\pi}$ is the group $\Phi_{B_\pi}$ where $B_{\pi}$ is an abelian variety in the corresponding isogeny class. So $\Phi_{\pi}$ is the $\absq$-stable subgroup of $\qbar^*$ generated by the conjugates of $\pi$ over $\bq$. We introduce the following notion that is weaker than that of isogeny between two abelian varieties over finite fields.

\begin{Def} \normalfont Following [Kow03], we say two abelian varieties $B_1$, $B_2$ over a finite field are \textit{isokummerian} if there exists some integer $N$ such that $k(B_1[d]) = k(B_2[d])$ for all $d$ prime to $N$.\end{Def} 
\vspace{-0.2cm}
The following lemma establishes an equivalence between this notion and that of the the multiplicative groups generated by the eigenvalues coinciding.

\begin{Prop} $\mathrm{([Kow03],\;Lemma\;8.1)}$ Two abelian varieties $B_1,B_2$ over a finite field are isokummerian if and only if $\Phi_{B_1} = \Phi_{B_2}$.\end{Prop}
\vspace{-0.2cm}
Clearly, two isogenous abelian varieties are isokummerian. The converse is true, for instance, when one of the two is an elliptic curve (on account of the lack of non-torsion units in imaginary quadratic fields). But it is not true in general even for ordinary abelian varieties, as evidenced by a class of counterexamples constructed in ([Kow06], Remark 2.5). However, for Weil numbers associated to ordinary abelian varieties, we have the following result that we will need.

\begin{Prop} $\mr{([Kow06],\;Proposition\;2.1)}$ Let $B$ be an abelian variety over a finite field $k$ of size $q$. Suppose the corresponding Weil number $\pi$ is such that $(\pi,\ov{\pi}) = (1)$ in the Galois closure $\bq(\W_B)$ of $\bq(\pi)$ over $\bq$. Suppose, further, that the exponent $2$ abelian Galois group $\Gal(\bq(\W_{B})/\bq(\W_{B})^+)$ is of rank $\dim B$. Then the only Weil $q$-integers in the multiplicative group $\Phi_B$ are the conjugates of $\pi$.\end{Prop}
\vspace{-0.2cm}
The next two lemmas will be useful for detecting the places of ordinary reduction. As shown in [Tat66], if a simple abelian variety $B_{\pi}$ is ordinary, the endomorphism algebra is commutative. While the converse is certainly false, we have the following two results in the reverse direction: 

\begin{Lem} Let $B_{\pi}$ be an absolutely simple abelian variety over a finite field $k$. Suppose:\\
$\mr{(i)}$ $\End^0_{\bark}(B_{\pi})$ is commutative\\
$\mr{(ii)}$ $p$ splits completely in $\bq(\pi)$.\\
Then $B_{\pi}$ is ordinary.\end{Lem}
\begin{prf} The Newton slopes of $B_{\pi}$ are given by $\frac{w(\pi)}{w(q)}$, where $w$ runs through all the places of $\bq(\pi)$ lying over $p$. Since $\End^0(B_{\pi})$ is commutative, it follows from Honda-Tate theory that the local invariants $[\bq(\pi)_w:\bq_p]\frac{w(\pi)}{w(q)}$ are integers. But since $p$ splits completely in $\bq(\pi)/\bq$ and hence, in the Galois closure $\bq(\W_{B_{\pi}})$, any prime of $\bq(\W_{B_{\pi}})$ lying over $p$ has local degree $1$ over $\bq$. So the Newton slopes $\frac{w(\pi)}{w(q)}$ are integers and hence, $B_{\pi}$ is ordinary.\end{prf}

\vspace{0.03cm}
\begin{Lem} Let $B_{\pi}$ be an absolutely simple abelian variety over a finite field $k$. Suppose:\\
$\mr{(i)}$ $\End^0_{\bark}(B_{\pi})$ is commutative\\
$\mr{(ii)}$ The Galois group $\Gal(\bq(\W_{B_{\pi}})/\bq)$ is such that the decomposition group of a prime of $\bq(\W_{B_{\pi}}$ lying over $p$ is a normal subgroup of $\Gal(\bq(\W_{B_{\pi}})/\bq)$.\\
Then $B_{\pi}$ is ordinary.\end{Lem}
\vspace{-0.15cm}
\noindent\underline{Remark} In particular, the second condition is fulfilled when the center $\bq(\pi)$ of $\End^0(B_{\pi})$ is an abelian extension over $\bq$. 

\begin{prf} Let $\mf{p}$ be a prime of $\bq(\W_{B_{\pi}})$ lying over $p$. If $\mf{p}$ has local degree $1$, then clearly all the Newton slopes are integral and hence, must be either $0$ or $1$. So we may assume without loss of generality that $p$ does not split completely in $\bq(\W_{B_{\pi}})$ and hence, the decomposition group of any prime lying over $p$ is non-trivial.

Choose a non-trivial element $\tau$ in the decomposition group $D_{\mf{p}}$. Since $\pi$ and its conjugates generate $\bq(\W_{B_{\pi}})$, $\tau$ cannot fix every conjugate of $\pi$. So we may assume without loss of generality that $\tau(\pi)\neq \pi$. Furthermore, the hypothesis implies that $D_{\mf{p}}$ is normal. So $\tau$ lies in the decomposition group of every prime of $\bqW$ lying over $p$ and hence, $\tau$ preserves the ideal $(\pi)$ in $\bq(\W_{B_{\pi}})$. So the algebraic number $\frac{\tau(\pi)}{\pi}$ has the valuation $0$ under every non-archimedean valuation and hence, must be a root of unity. Write $\tau(\pi) = \zeta_{_N}\pi$. Then the endomorphism algebra ${\End^0(B\times_{\bF _{q}} \bF _{q^N})}$ is non-commutative since the characteristic polynomial of $B\times_{\bF _{q}} \bF _{q^N}$ has the zero $\pi^N$ with multiplicity $\geq 2$, a contradiction. So $p$ splits completely in $\bq(\W_{B_{\pi}})/\bq$ and by the preceding lemma, it follows that $B_{\pi}$ is ordinary. \end{prf}

\vspace{-0.4cm}
\subsection{\fontsize{11}{11}\selectfont Newton Polygons}
\vspace{-0.3cm}
Let $A$ be an abelian variety weakly of Mumford type of rank $6$. Let $v$ be a place of good reduction, so that the Newton polygon of $A_v[p]$ lies above the Hodge polygon. With techniques virtually identical to those of ([No00], Section 3), we may show that the Newton slopes are given by\vspace{-0.15cm}$$\frac{1}{2}\pm a_1\pm a_2\pm a_3\pm a_4\pm a_5$$ where the $a_i$ are rational numbers between $0$ and $1$. Using the facts that the Newton slopes lie between $0$ and $1$, have integral break points and the Newton polygons of abelian varieties are admissible, we can then deduce that the slope $\frac{1}{2}$ occurs with multiplicity $0, 8, 12$ or $16$. This will be useful in the proof of Proposition 4.8.
\vspace{-0.4cm}
\subsection{\fontsize{11}{11}\selectfont Reductions of abelian varieties weakly of Mumford type}
\vspace{-0.3cm}
Let $A$ be a $16$-dimensional (equivalently, rank $6$) abelian variety weakly of Mumford type over a number field $F$. Enlarging $F$ if necessary, we may assume without loss of generality that $A$ has good reduction everywhere. Let $v$ be a non-archimedean place with $p=\mathrm{char}(v)$. We choose a prime $l\neq p$ such that the derived subgroup $G_l^{\der}$ is $\bq_l$-simple and $(G_l,V_l)$ is of Mumford type. As seen in lemma 3.4, the set of primes fulfilling the first condition has positive density and since $\mb{S}_A$ has density one, the intersection has positive density as well. In particular, there are infinitely such primes.

\begin{Prop} With the setup as above, the reduction $A_{\ov{k}}$ has one of the following simple decompositions up to isogeny:

\noindent $1.$ $(A^{(1)})^{16}$ where $A^{(1)}$ is a supersingular elliptic curve\\
$2.$ $A^{(1)}\times A_1^{(5)}\times A_2^{(5)}\times A_2^{(5)}$ with Newton polygon $0, 5\times \frac{1}{5}, 10\times \frac{2}{5}, 10\times\frac{3}{5}, 5\times \frac{4}{5}, 1$\\
$3.$ $A^{(1)}\times A^{(5)}\times A^{(10)}$ with $A^{(1)}$ supersingular\\
$4.$ $A^{(1)}\times A^{(5)}\times A^{(10)}$ with $A^{(1)}$ ordinary\\
$5.$ $A^{(1)}\times A_1^{(5)}\times A_2^{(5)}\times A_3^{(5)}$ with each component ordinary and the last two isokummerian\\
$6.$ Simple

where the superscript is the dimension of the simple abelian variety and distinct subscripts indicate non-isogeny. Furthermore, in the last three cases, $\End(A_{\bark})$ is commutative. If the base field is enlarged so that $F=F_A^{\conn}$, the last possibility occurs with Dirichlet density one.	\end{Prop}
\vspace{-0.2cm}
\noindent\underline{Remark} In the course of this proof, we will also  determine the possibilities for the endomorphism algebra $\End^0(A_{\ov{k}})$ and state them in Proposition 4.11. That is a more concrete description of the simple decomposition since the dimension of a simple abelian variety can be easily retrieved from the endomorphism algebra via Honda-Tate theory. Furthermore, we also look into the Galois group $\GalW$ which provides some insight into the reduction $A_v$.
\vspace{0.1cm}
\begin{prf} Let $$\rho_{v,l}:\absfv\cong D_v\hra \absf \lra G_l(\bz_l)$$ be the local $l$-adic Galois representation attached to $A$. Replacing $F$ by a finite extension if necessary, we may assume the following:

\noindent 1. $q:=|k|$ is a perfect square.\\
2. All the $l$-adic monodromy groups are connected.\\
3. All elements of $\rho_{v,l}(\absfv)$ are congruent to the identity mod $l^2$.\\
4. All simple factors of $A_{\ov{k}}$ are defined over $k$.

In particular, the simple components of $A_v$ are absolutely simple. So for any Weil number $\pi$ corresponding to a simple component of $A_v$, the inclusion $\bq(\pi^d) \sub \bq(\pi)$ is an equality for any integer $d$. The second condition implies that any element of $\rho_{v,l}(\absfv)$ lifts uniquely to an element of $\ti{G}$. The third ensures that the group generated by the eigenvalues of the Frobenius does not contain any roots of unity other than $1$.

Now, there exists a universal simple cover $\Phi:\tilde{G_l}\lra G_l$ with a central isogeny $\Phi$. From the third assumption above, $\rho_{v,l}$ lifts uniquely to a representation $\tilde{\rho}_{v,l}:\absfv\lra \ti{G_l}(\bz_l)$ with the same congruence condition. Because of the uniqueness of the lift, $\ti{\rho}_{v,l}$ is a group homomorphism. Let $\mr{Fr_v}$ be a geometric Frobenius. Set $\pi = \ti{\rho}_{v,l}(\mr{Fr_v})\in \ti{G_l}(\bq_l)$ and define $\wti{T}\sub\ti{G_l}$ as the Zariski closure of the subgroup of $\ti{G_l}(\bq_l)$ generated by $\pi$. The congruence condition on $\ti{\rho}$ ensures that the subgroup of $\qbar^*$ generated by the eigenvalues of $\pi$ does not contain any roots of unity other than $1$ and hence, $\wti{T}$ is connected. So $\wti{T}$ is a maximal torus in $\ti{G_l}$.

Now $\sqrt{q}\in \bz$ by construction. Set $\al:= \frac{\pi}{\sqrt{q}}\in \tilde{G}_l(\bq_l)$ and let $\wti{T'}$ be the Zariski closure of the subgroup of $\tilde{G}_l$ generated by $\alpha$. Then $\wti{T}\cong \bg_{m,\bq_l}\times \wti{T'}$ with $\wti{T'}$ a maximal torus in $\ti{G}^{\der}$. We choose a maximal torus $\ov{T}_{\bq_l}$ of $\ti{G_l}^{\der}$ (over $\bq_l$) containing $\wti{T'}$. Since $(G_l,V_l)$ is of Mumford type of rank $6$, $\ov{T}$ is of rank $5$. Let $\Psi$ be the associated root system of $\wti{T'}$. Let $\mr{Out}(\Psi) = \mr{Aut}(\Psi)/W(\Psi)$ be the group of outer automorphism group of $\Psi$, which is isomorphic to $S_5$. There exists an isomorphism $X^*(\ov{T})\cong \bz^{5}$ such that the weights of the representation of $\ov{T}$ on $V_l$ correspond to the elements ${(\pm 1,\cdots,\pm 1)\in \bz^{5}}$. The action of $\absql$ on $X^*(\ov{T})$ induces an action on $\bz^{5}$. This action factors through the group $\{\pm 1\}^5\rtimes S_{5}$ and contains a $5$-cycle, which we may assume without loss of generality to be the cycle $(1\;2\;3\;4\;5)$.

Evaluation at the element $\al\in \ov{T}(\bq_l)$ defines a $\absql$-equivariant map \vspace{-0.1cm}$$\ev :X^*(\ov{T})\lra \qbar_l^*,\;\;\;\chi\mapsto \chi(\alpha).$$ But the characteristic polynomial of the Frobenius lies in $\bz[X]$ and is independent of $l$. Hence, the eigenvalues of $\rho_{v,l}(\Fr_v)$ lie in $\qbar$ where $\qbar$ is some algebraic closure of $\bq$ embedded in $\qbar_l$. Since the weights of the representation of $\ov{T}$ on $V_l$ correspond to the elements $(\pm 1,\pm 1,\pm 1,\pm 1,\pm 1)\in \bz^{5}$, the images $\ev(\pm 1,\cdots, \pm 1)$ are the eigenvalues of $\alpha$ on $V_l$. It follows that $\ev(\pm 1,\cdots,\pm 1)$ are in $\qbar\sub\qbar_l$. So the image of $\ev$ lies in $\qbar^*$. As $\wti{T'}$ is the Zariski closure of the subgroup of $\ti{G_l}(\bq_l)$ generated by $\alpha$, it follows that $\ker(\ev)$ is the kernel of the natural surjection  $X^*(\ov{T})\lra X^*(\wti{T'})$. So we may identify $X^*(\wti{T'})$ with $\bz^{5}/\ker(\ev)$.

For any element $(x_{_1},\cdots,x_{5})$, denote its image in $X^*(\wti{T'})$ by $(x_{_1},\cdots,x_{5})'$. The map $\ev$ induces an action of $\absq$ on $X^*(\wti{T'})$ extending the action of $\absql$. This action stabilizes the elements $\{(\pm 1,\cdots, \pm 1)\}\sub X^*(\wti{T'})$. Since all the eigenvalues of $\alpha$ have absolute value $1$, it follows that $\ev(P)\ov{\ev(P)}=1$ for each vertex $P$ of the hypercube $\{(\pm 1,\cdots, \pm 1)\}$. So the complex conjugation acts by multiplication by $-1$ on $X(\wti{T'})$.

We fix an embedding $\qbar\hra \qbar_p$ which induces a $p$-adic valuation $\mathfrak{v}$ on $\qbar$ normalized by $\mathfrak{v}(q) = 1$  and define $\phi = \mathfrak{v}\circ \mr{ev'}:X^*(\wti{T'})\lra \bq$. This map is $\bz$-linear but not necessarily $\absq$-equivariant or injective. The Newton slopes are then given by $\frac{1}{2}+\phi(\pm 1,\pm 1,\pm 1,\pm 1,\pm 1)$. Although $\phi$ might have a non-trivial kernel, the following lemma shows that it does not have a $\absq$-stable subset. 

\begin{Lem} Let $x\in X^*(\wti{T'})$. If $\phi(\tau(x)) =0$ for all $\tau\in \absq$, then $x=0$.\end{Lem}

\begin{prf} Since the subgroup of $X^*(\wti{T'})$ generated by $\{(\pm 1,\cdots,\pm 1)\}$ is of finite index, we may assume without loss of generality that $x$ lies in the subgroup. As $\tau$ runs through $\absq$, $\phi(\tau(x))$ runs through all $p$-adic valuations of $\ev'(x)$. So these are all $0$ and hence, all absolute values of $\ev'(x)$	are $0$. So $\ev'(x)$ is a root of unity, which implies that $\ev'(x)=1$. Since $\ev'$ is injective, this is a contradiction. \end{prf}

First, note that $\bz^5/\ker(\ev)$ is torsion-free since the image of $\ev$ contains no roots of unity other than $1$. Now $\ker(\ev)$ is stable under the $\absql$-action and the image of $\absql$ in $\{\pm 1\}\rtimes S_{5}$ contains a $5$-cycle. So $\ker(\ev)$ is a $\absql$-stable sublattice of $\bz^5$ with the quotient $\bz^5/\ker(\ev)$ torsion-free. So $\ker(\ev)$ is a cyclic sublattice with a torsion-free quotient and hence, must be one of the following:

\noindent 1. $\bz^{5}$

\noindent 2. $\{(x_{_1},\cdots,x_{5}):x_{_1}+\cdots+x_5 =0\}$

\noindent 3. $\{(x,\cdots,x): x\in \bz\}$

\noindent 4. $ 0$

\noindent(see the appendix on cyclic lattices for an explanation. Note that in the case of good ordinary reduction - which occurs with density one- only case $4$ can occur.)

In \textbf{Case 1}, all eigenvalues of $\Fr_v$ are $\sqrt{q}$. So all Newton slopes are $\frac{1}{2}$. Hence, $A$ is isogenous to the $16$-th power of the supersingular elliptic curve and $\End^0(A_{\ov{k}}) = \Mat_{16}(\bq_{p,\infty})$, where $\bq_{p,\infty}$ is the quaternion algebra over $\bq$ ramified only at $p$ and $\infty$. Clearly, $\GalW = \{1\}$.

In \textbf{Case 2}, $X^*(\wti{T'})$ is a free $\bz$-module of rank $1$. Since none of the weights $(\pm 1,\pm 1,\pm 1,\pm 1,\pm 1)$ are in the kernel of $\ev$, none of the simple components correspond to $\sqrt{q}$.

Write $a_1 = (1,1,1,1,1)', a_2 = (-1,1,1,1,1)', \cdots, a_6  = (-1,-1,-1,-1,-1)'\in X^*(\wti{T'})$. Then \vspace{-0.15cm}$$a_1-a_2 = a_2-a_3 = a_3-a_4 = a_4-a_5 = a_5-a_6$$ and $\phi(a_1) = -\phi(a_6)$. Since the only conjugates of $\ev(a_1)$ are $\ev(a_1)$ and $\ov{\ev(a_1)}$, lemma 4.9 implies that $\phi(a_1)\neq 0$. Hence, $\phi(a_1) = 5\lambda$, $\phi(a_2) = 3\lambda$, $\phi(a_3) = \lambda$, $\phi(a_4) = -\lambda$, $\phi(a_5) = -3\lambda$, $\phi(a_6) = -5\lambda$ for some $\lambda$. So the Newton slopes are $\frac{1}{2}\pm 5\lambda$, $5\times \frac{1}{2}\pm 3\lambda$, $10\times \frac{1}{2}\pm \lambda$. Since the slopes lie between $0$ and $1$ and the Newton polygon has integral break points, it follows that $\lambda = \frac{1}{10}$. Hence, the 	Newton polygon of $A_v$ is $0, 5\times \frac{1}{5}, 10\times \frac{2}{5}, 10\times\frac{3}{5}, 5\times \frac{4}{5}, 1$.

Now $\rho_{_l}(\Fr_v)$ has six distinct $\ev(a_1),\cdots,\ev(a_6)$ eigenvalues on $V_l$.  We have $\ev(a_1) = \overline{\ev(a_6)}$ is of degree $2$ and hence, $\ev(a_1)\sqrt{q}$ is a Weil $q$-integer of degree $2$. So it corresponds to an ordinary elliptic curve $A^{(1)}$. Write $A_v = A^{(1)}\times A'$. Since $A'$ has Newton polygon $5\times \frac{1}{5},10\times \frac{2}{5}, 10\times \frac{3}{5}, 5\times \frac{4}{5}$, every simple component has dimension divisible by $5$. Since $\GalW = \{\pm 1\}\times H$ for some $H\sub S_5$, it follows that  
$A'$ has a simple component of dimension $5$. So $\absq$ has precisely three orbits in $\W_{A_v}$ and hence, $A'$ has two simple components (possibly with multiplicities greater than one).

Now suppose, by way of contradiction, that $A'$ has a simple component $B$ of dimension $10$. Its Weil number is of degree $2$ over $\bq$ and hence, $20 = 2[\End^0(B):\Cent(\End^0(B))]^{1/2}$. So $[\End^0(B):\Cent(\End^0(B))]^{1/2} = 10$. But the exact sequence \vspace{-0.1cm}$$0\lra \Br(L|K)\xra{\;\;i\;\;} \bigoplus_w \Br(L^w|K_w)\xra{\al\mapsto \sum_v \inv_v(\al)} \frac{1}{\lcm(\{n_w\}_w)}\bz/\bz\lra 0$$ \big(where $L^w$ is the completion of $L$ at any prime lying over $w$ and $n_w = [L^w:K_w]$\big) for Galois extensions $L/K$ of global fields implies that $[\End^0(B):\Cent(\End^0(B))]^{1/2}$ is the least common denominator of the local invariants $\inv_w(\End^0(B))$. But this integer, in turn, divides the least common denominator of the Newton slopes, which is $5$. This yields a contradiction.

Thus, $A'$ has two simple components $B_1$ and $B_2$ each of dimension $5$ and of multiplicities $1$ and $2$ respectively. Furthermore, both $\End^0(B_1)$ and $\End^0(B_1)$ are $25$ dimensional central division algebras over imaginary quadratic fields in which $p$ splits.

In \textbf{Case 3}, $X^*(\wti{T'})$ is a free $\bz$-module of rank $4$. Now $(1,1,1,1,1)' = (-1,-1,-1,-1,-1)' = 0$. 	The Weil $q$-integer $\sqrt{q} = \ev(1,1,1,1,1)\sqrt{q}$ corresponds to a supersingular elliptic curve $A^{(1)}$ over $k$. Furthermore, the other $30$ eigenvalues $\W_{A_v}\setminus \ev(\pm(1,1,1,1,1))$ are all distinct. So the endomorphism algebra of the $15$ dimensional complement $A'$ of $A^{(1)}$ is a CM algebra of dimension $30$. The image of $\absq$ in $\{\pm 1\}\rtimes S_5$ contains a cycle of length $5$ which induces a cyclic permutation of $\{\ev(-1,1,1,1,1),\cdots,\ev(1,1,1,1,-1) \}$. Since complex conjugation acts by inversion, $\absq$ acts transitively on $\{\ev(\pm(-1,1,1,1,1)),\cdots,\ev(\pm(1,1,1,1,-1)) \}$. Similarly, $\absq$ acts transitively on \vspace{-0.15cm}$$\{\ev(\pm(-1,-1,1,1,1)),\cdots,\ev(\pm(-1,1,1,1,-1)) \}$$ and $$\{\ev(\pm(-1,1,-1,1,1)),\cdots,\ev(\pm(1,-1,1,1,-1)) \}.$$

Set $\pi_{_1} = \ev(-1,1,1,1,1)\sqrt{q}$, $\pi_{_2} = \ev(-1,-1,1,1,1)\sqrt{q}$, $\pi_{_3} = \ev(-1,1,-1,1,1)\sqrt{q}$. Each $\pi_{_i}$ ($i=1,2,3$) is a Weil $q$-integer of degree a multiple of $10$. Let $A_i$ be the corresponding simple abelian variety over $\bF _q$ (up to isogeny), $L_i$ the field $\bq(\pi_{_i})$ and $\al_i = \frac{\;\pi_{_i}}{\sqrt{q}}$. Since $\End^0(A_i) = L_i$, it follows that the local invariants $[L_{i,v}:\bq_p]v(\tau(\pi_{_i}))\in\bz$ for all $\tau\in\absq$.

Furthermore, we have $\Phi_{\pi_{_1}} = \Phi_{\pi_{_2}} = \Phi_{\pi_{_3}}$ and $\Phi_{\pi_{_0}} = \la \sqrt{q} \ra\sub \qbar^* .$ Now, $\GalWr$ is a transitive subgroup of $S_5$. So it must be one of the following groups: \vspace{-0.15cm}$$\{\la(1\;2\;3\;4\;5)\ra,\la (1\;2\;3\;4\;5), (1\;5)(2\;4)\ra,\la (1\;2\;3\;4\;5), (1\;2\;5\;4)  \ra, A_5,S_5 \}.$$ So every simple component of $A'$ has dimension a multiple of $5$. Suppose by way of contradiction that every simple component is a five-fold. Then, in particular, $\pi_{_2}$ is not conjugate to $\pi_{_3}$. Hence, $\GalWr$ is cyclic of order $5$ and $\GalW$ is cyclic of order $10$. Then $\End^0(A')$ is commutative and it follows from lemma 4.6 that every simple component of $A'$ is ordinary. Hence, the Newton polygon is $2\times \frac{1}{2}, 15\times 0, 15\times 1$ which is not possible as noted in subsection 4.3. So $\GalWr$ is one of the other four possibilities and $\pi_{_2}=\ev(-1,-1,1,1,1)\sqrt{q}$ is conjugate to $\pi_{_3} = \ev(-1,1,-1,1,1)\sqrt{q}$. Hence, the simple component corresponding to $\pi_{_2}$ is of dimension $10$.

The endomorphism algebra $\End^0(A_v)\cong \bq_{p,\infty}\times F_{10}\times F_{20}$ where the subscript is the degree of the CM field and $F_{10}\hra F_{20}$.

In \textbf{Case 4}, there are $32$ distinct eigenvalues and the endomorphism algebra is a CM algebra of degree $32$. An argument identical to the one in case 3 shows that $\absq$ acts transitively on the sets $\{\ev(\pm(1,1,1,1,1))\}$, ${\{\ev(\pm(-1,1,1,1,1)),\cdots,\ev(\pm(1,1,1,1,-1)) \}}$,\\ ${\{\ev(\pm(-1,-1,1,1,1)),\cdots,\ev(\pm(-1,1,1,1,-1)) \}}$ and ${\{\ev(\pm(-1,1,-1,1,1)),\cdots,\ev(\pm(1,-1,1,1,-1)) \}}$. So the simple component corresponding to the Weil number $\pi_{_0}$ has dimension $\equiv 1\pmod 5$ and the other simple components have dimensions divisible by $5$. Furthermore, we have \\$\GalW\sub \{\pm 1\}^j\rtimes H$ and since $H$ fixes $\pi_{_0} = \ev(1,1,1,1,1)\sqrt{q}$, we have $\bq(\pi_{_0})\sub \bqW^H$. So $[\bq(\pi_{_0}):\bq] = 2^i$ for some $1\leq j\leq 5$ and $2^{i-1}\equiv 1\pmod 5$, which implies that $i=1$ or $5$.

Also, we have \vspace{-0.15cm}$$\Phi_{\pi_{_0}}\sub \Phi_{\pi_{_1}} \sub \Phi_{\pi_{_2}} = \Phi_{\pi_{_3}} = \Phi_{A_v}.$$ Furthermore, we have the inclusion $\bq(\pi_0)^{\Gal}\sub \bq(\pi_1)^{\Gal} = \bq(\pi_2)^{\Gal} = \bq(\pi_3)^{\Gal}$ of Galois closures and $[\bq(\pi_0):\bq]\geq 2$, $[\bq(\pi_1):\bq]\geq 10$, $[\bq(\pi_2):\bq]\geq 10$, $[\bq(\pi_3):\bq]\geq 10$. Thus, the special fiber $A_v$ has one of the following simple decompositions in \textbf{Case 4}:

\textbf{(i)} Simple of dimension $16$\\
In this case, we have $[\bq(\pi_{_0}):\bq]=32$. Write $\GalW = \{\pm 1\}^5\rtimes H$. Since $H$ fixes $\pi_{_0}$, it follows that $A_v$ has complex multiplication by the field $\End^0(A_v) \cong \bq(\pi_{_0})$, the subfield of $\GalW$ fixed by $H$.

In all other cases, $\pi_{_0}$ corresponds to an ordinary elliptic curve. Suppose, by way of contradiction, that $A_v$ is not simple and $\Gal(\bqW/\bqWr)$ is of rank $5$. Choose an element $\tau_j\in \Gal(\bqW/\bqWr)$ that fixes $\si^j(\pi_{_2})$ and swaps $\si^i(\pi_{_2})$ and $\si^j(\ovp_2)$ for $i\neq j$, ($0\leq i,j\leq 4$). We have \vspace{-0.1cm}$$\pi_{_0} = \pi_{_2}\si(\ovp_2)^{-1}\si^2(\pi_{_2})\si^3(\ovp_2)^{-1}\si^4(\pi_{_2}).$$ So $\pi_{_0}\pi_0^{\tau_{_0}} = \pi_{_2}^2$ and hence, $\pi_{_0}\pi_{_0}^{\tau_{_0}} = q$ or $\pi_{_0}^2$ depending on whether $\tau_{_0}$ fixes $\pi_{_0}$. So either $\pi_{_2}^2 = q$ or $\pi_{_2}^2 = \pi_{_0}^2$, neither of which is possible in \textbf{Case 4}. So $\Gal(\bqW/\bqWr)$ is of order $2$.

Let $A'$ be the $15$-dimensional complement of $A_{\pi_{_0}}$ in $A_v$ and let $B$ be an odd dimensional simple component of $A'$. Let $\pi_{_i}$ be the corresponding Weil number (where $i=1$, $2$ or $3$). Then $\bq(\pi_{_i}+\ov{\pi}_{_i})$ is of odd degree over $\bq$ and hence, lies in the fixed field of $\{\pm 1\}$. Since $\{\pm 1\}$ is a normal subgroup of $\GalW$, the Galois closure $\bqWr$ of $\bq(\pi_{_i}+\ov{\pi}_{_i})$ also lies in the fixed field of $\{\pm 1\}$. Hence, $\GalW = \{\pm 1\}\times \GalWr$.

\textbf{(ii)} $A_{v} = A^{(1)}\times A^{(5)}\times A^{(10)}$\\
In this case, $\pi_2,\pi_3$ lie in the same $\absq$-conjugacy class. So $\GalW = \{\pm 1\}\times H$ where $H$ is a transitive subgroup of $S_5$ of order divisible by $10$.

\textbf{(iii)} $A_{v} = A^{(1)}\times A_1^{(5)}\times A_2^{(5)}\times A_3^{(5)}$\\
In this case, $\bq(\pi_{_0})$ is imaginary quadratic with $p$ split in $\bq(\pi_{_0})/\bq$ and the elements $\pi_{_i}$ ($i=1,2,3$) are pairwise non-conjugate and of degree $10$. Since $\pi_{_2}$, $\pi_{_3}$ are non-conjugate, it follows that $\GalWr$ is cyclic of order $5$ and $\GalW$ is cyclic of order $10$. So $p$ splits completely in $\bqW/\bq$ and hence, by lemma 4.7, every simple component is ordinary.

\textbf{(iv)'} $A_v = A^{(1)}\times A^{(15)}$\\
In this case, $\pi_{_1},\pi_{_2},\pi_{_3}$ lie in the same conjugacy class and $\pi_{_0}$ is of degree $2$. Furthermore, $30$ divides $|\GalW|$ and hence, $\GalW =\{\pm 1\}^j\rtimes H$ with $j = 1$ or $5$ and $H = A_5$ or $S_5$. Since the set $\{\si^i(\pi_{_1}),\si^i(\ov{\pi}_{_1}):1\leq i\leq 5\}$ is stable under the action of $H$ and the complex conjugation, it follows that $j=5$, a contradiction.\end{prf}

We briefly describe how the techniques in the preceding theorem may be used for higher dimensions. Consider an abelian variety $A$ weakly of Mumford type of rank $N+1$ with $N$ a prime. As explained in appendix 6.2, $\ker(\ev)$ has the same four possibilities as above. 

Clearly, if $\ker(\ev)$ is the full $\bz^N$, then all the eigenvalues are $\frac{1}{2}$ and the reduction is isogenous to the $2^{N-1}$-th power of the supersingular elliptic curve. 

If $\ker(\ev)$ is the index $1$ sublattice $\{(x_{_1},x_{_2},\cdots,x_{_N}):\sum\limits_{i=1}^N x_i=0\}$, then $A_v$ has precisely $N+1$ distinct eigenvalues. Furthermore, with techniques identical to those in the preceding proposition, we may show that the Newton polygon of $A_v$ is \vspace{-0.1cm}$$ 2^{N-j}{N\choose j}\times \frac{j}{N},\;\;\;\;0\leq j\leq N.$$ So every simple component of $A_v$ other than the ordinary elliptic curve has dimension a multiple of $N$. Furthermore, for any such simple component $B$, the integer $[\End^0(B):\Cent(\End^0(B))]^{1/2}$ is the least common multiple of the Hasse invariants $\inv_w([\End^0(B)])$ and hence, must divide $N$. Since $N$ is a prime by assumption, it follows that $[\End^0(B):\Cent(\End^0(B))] =N^2$. So each $\End^0(B)$ is a $N^2$-dimensional central division algebra over an imaginary quadratic field in which $p$ splits. Thus, by Honda-Tate theory, every simple component other than the ordinary elliptic curve has dimension $N$ and $\GalW$ is abelian of exponent $2$. The simple component of $A_v$ corresponding to the Weil number $\ev(\underbrace{-1,\cdots,-1}_{j},1,\cdots,1)\sqrt{q}$ occurs with multiplicity $\frac{1}{N}{{N}\choose {j}}$ and has precisely two Newton slopes $\frac{j}{N}, \frac{N-j}{N}$, each with multiplicity $N$.  

If $\ker(\ev) = \{(x,x,\cdots,x):x\in \bz\}$, then $\GalW = \{\pm 1\}\times H$ where $H$ is some subgroup of $S_{_N}$ containing a $N$-cycle. Clearly, $H$ determines the $\absq$-orbits of $\WA\setm \{\sqrt{q}\}$ and hence, determines the simple decomposition of the complement of the supersingular elliptic curve in $A_v$.

We conclude this subsection by stating a result for the case when the map $\ev$ is surjective. As before, we denote the Weil number $\ev(1,\cdots,1)\sqrt{q}$ by $\pi_{_0}$.

\begin{Prop} Let $N$ be a prime. Let $A$ be an abelian variety over a number field, weakly of Mumford type of rank $N+1$ and let $v$ be a non-archimedean place of $F$. Then $[\bq(\pi_{_0}):\bq] = 2^i$ for some $1\leq i\leq N$ such that $2^{i-1}\equiv 1\pmod {N}$ and $A_v$ is simple if and only if $[\bq(\pi_{_0}):\bq] = 2^{N}$.\end{Prop}
\begin{prf} By lemma 3.5, we may enlarge the field if necessary to ensure good reduction at $v$. Furthermore, the kernel of the map $\ev:X^*(T)\lra \Phi_{A_v}$ has the four possibilities listed in appendix 6.2. Unless $\ker(\ev) = 0$, we see that $\pi_{_0}$ corresponds to an elliptic curve. So we may assume without loss of generality that $\ker(\ev)=0$. Hence, $P_{A_v}(X)$ is separable and $\End^0(A_v)$ is a CM algebra of degree $2^{N}$.

Now, note that $\GalW = \{\pm 1\}^{j}\rtimes H$ where some $1\leq j\leq N$ and $H$ is a transitive subgroup of $S_{_N}$. Since $\pi_{_0}$ is fixed by $H$, we have $\bq(\pi_{_0})\sub \bqW^H$ and since\\ $[\bqW^H:\bq] = 2^j$, it follows that {$[\bqW^H:\bq] = 2^{i-1}$} for some $i\leq j$. Furthermore, since $N$ is assumed to be a prime, $H$ contains a $N$-cycle. Hence, every simple component of $A_v$ corresponding to a Weil number non-conjugate to  $\pi_{_0}$ has dimension a multiple of $N$ and the component corresponding to $\pi_{_0}$ has dimension $\equiv 1 \pmod{N}$. Thus, $2^{i-1}\equiv 1\pmod {N}$.

Since the eigenvalues are distinct, it is clear that $A_v$ is simple if and only if $[\bq(\pi_{_0}):\bq] = 2\;\dim A_v = 2^{N}$.\end{prf}

In particular, we note that at places where the abelian variety has simple reduction, the eigenvalues of the Frobenius are all distinct and the endomorphism ring is an order in a CM field. So the image of $\Fr_v$ in $\End(T_l(A))$ generates a maximal torus of $G_l$.
\vspace{-0.3cm}
\subsection{\fontsize{11}{11}\selectfont Endomorphism Algebras}
\vspace{-0.3cm}
We now list the possibilities for $\End^0(A_{\ov{k}})$ that we derived in the course of the proof of Theorem 4.8. Given the division algebra $\End^0(Y)$ with center $C$ for a simple component $Y$, we can easily retrieve the dimension through the equality $\dim Y = \frac{1}{2}[\End^0(Y):C][C:\bq]$. Hence, this is a more concrete characterization than the dimensions of the simple components, even though the endomorphism algebras do not determine the precise isogeny classes for abelian varieties of genus greater than one.

We introduce some notations for brevity. The center of the endomorphism algebra of a simple component of $A_{\ov{k}}$ is some field $\bq(\sqrt{q}\alpha)$ with $\sqrt{q}\alpha$ a Weil $q$-integer. So it is either $\bq$ (which only happens in the supersingular case) or a CM field. We use the symbol $\bq_{p,\infty}$ for the unique quaternion algebra over $\bq$ ramified only at $p$ and $\infty$. We will use $F_{2i}$ to mean a CM field of degree $2i$ and $D^{2i,j}$ to mean the central division algebra of dimension $j^2$ over a CM field $E_{2i}$ of degree $2i$ such that:\\
- the rational prime $p$ splits completely in $E_{2i}$;\\
- $D^{2i,j}$ has Hasse invariants\vspace{-0.1cm} $$\inv_w(D^{2i,j}) = \begin{cases} 0 & \text{ if }w\nmid p\\
\frac{w(\pi)}{w(q)} [E_{2i,w}:\bq_p]\;\pmod{1} & \text{ if } w|p.\end{cases}$$

\vspace{-0.2cm}

\begin{Prop} For an abelian variety $A$ weakly of Mumford type of rank $6$ over a number field $F$, $\End^0(A_{\ov{k}})$ has the following possibilities:

\noindent $1.$ $\Mat_{16}(\bq_{p,\infty})$ $(\text{Supersingular case})$\\
$2.$ $F_2\times D_1^{2,5}\times \Mat_2(D_2^{2,5})$\\ 
$3.$ $\bq_{p,\infty}\times F_{10}\times F_{20}$\\
$4.$ $F_2\times F_{10}\times F_{20}$ \\
$5.$ $F_2\times F_{10}\times F_{10}\times F_{10}$ with $\Gal(F_{10}/\bq) = \bz/10\bz$, $F_2\sub F_{10}$ and $p\in \mc{S}(F_{10}|\bq)$\\
$6.$ $F_{32}$  

If the base field is enlarged so that $F=F_A^{\conn}$, the last possibility occurs with Dirichlet density one.\end{Prop}

\begin{prf} This is shown in the proof of Proposition 4.8. The density of places with simple ordinary reduction will be shown in section $5$.\end{prf}

\noindent\underline{Remark} We note that in the cases where the endomorphism algebra is a CM algebra, the equality of endomorphism algebras holds over the smallest extension of $k$ over which the reduction is well-defined.

\vspace{-0.3cm}	
\subsection{\fontsize{11}{11}\selectfont An existence result in the ordinary case}
\vspace{-0.3cm}
In the preceding subsection, we saw that for an abelian variety weakly of Mumford type of rank $6$, the reduction at a place of good ordinary reduction is either simple or of the form $A^{(1)}\times A'$, where $A^{(1)}$ is an elliptic curve. In this section, we show the existence of the latter possibility. Except for Proposition 4.13, the techniques are similar to those from [No01] and hence, we only provide a sketch.

\underline{The tori $T_L'$ and $T_L$.} We introduce some notation here. Let $L$ be a CM field and let $L^{^+}$ be its totally real subfield of index $2$. The norm map $N_{L/L^{^+}}:L^*\lra L^{^{+*}}$ defines a map\\ $\Res_{L/\bq}(\Gm)\lra \Res_{L^{^+}/\bq}(\Gm)$. We define $T_L'$ to be the kernel of the map and $T_L$ to be the pre-image of $\bg_{m,\bq}$.

Now let $K$ be a totally real field of degree $2n+1$, with $2n+1$ a prime. For simplicity, we assume $K$ is Galois over $\bq$. Let $D$ be a quaternion algebra over $K$ split at exactly one real place and with trivial corestriction in $\Br(\bq)$. By the Skolem-Noether theorem, we have $(D^*)^{\der} = \{d\in D^*:d\ov{d}=1\}$. The datum $((D^*)^{\der},h)$ gives a Shimura curve that is of Hodge type but not of PEL type. We refer the reader to ([SZ12], Section 2) for more details. 

\begin{Def}\noindent For a Shimura variety $(G,X)$ of Hodge type, we say a point $x$ is \textit{generic} if the Mumford-Tate group of $x$ is the same as that of $X$ and \textit{special} if the Mumford-Tate group is a torus.\end{Def}

Let $X$ be an abelian variety strongly of Mumford type that is a generic point on the Shimura curve given by $((D^*)^{\der},h)$. Then $X$ has a model over the field $K$ ([No01], [SZ12]). The following theorem of Tate guarantees the existence of an intermediate CM field between a totally real field and a quaternion algebra over it.

\begin{Prop} $\mr{(Tate)}$ Let $K$ be a totally real field and let $D$ be a $d^2$-dimensional central division algebra over $K$. If $d$ is even, there exists a CM field $L$ with an embedding in $D$ such that $[L:K] = d$.\end{Prop}
\vspace{-0.2cm}
A proof may be found in ([CCO13], Chapter 1). We will need the following variant for our purposes, a proof for which we provide in appendix 6.3.

\begin{Prop} Let $K$ be a totally real field Galois over $\bq$ and a $D$ a $4d^2$-dimensional central division algebra over $K$. Let $\wti{K}$ be the Galois closure over $\bq$. If $2\nmid d$ and $\gcd(2d,[\wti{K}:\bq]) =1$, there exists a totally imaginary number field $E$ cyclic of degree $2d$ over $\bq$ such that $[EK:K] = 2d$ and $EK$ has an embedding in $D$.\end{Prop} 
\vspace{-0.2cm}
In particular, there exists an imaginary quadratic field $E$ such that $L=EK$ has an embedding in the quaternion algebra $D$. The inclusion $L^*\sub D^*$ induces morphisms $\ti{\rho}':T'_L\lra \ti{G}'$ and ${\rho}' = N\circ \ti{\rho}':T'_L\lra G'$. Let $\phi_1,\cdots,\phi_{2n+1},\ov{\phi}_{1},\cdots,\ov{\phi}_{2n+1}$ be the embeddings $L\hra \bc$. Let $\phi_{2n+1}$ be the embedding corresponding to the archimedean place where $D$ is split. We write $[\phi_j]$, $[\ov{\phi_j}]$ for the induced characters of $L^*$, of $T_L^*$ and of $T_L'$ respectively. So $[\phi_1],\cdots,[\phi_{2n+1}]$ is a basis for $X^*(T_L')$. We fix an isomorphism $L\otimes_{\bq} \br\cong \bc^{2n+1}$ such that the composite of $L\hra L\otimes \br$ with the projection on the $i$-th factor $\bc$ induces the embedding $\phi_i:L\hra \br$. We deduce an isomorphism $L_{\br}^*\cong S^{2n+1}$. Define $h$ to be the composite \vspace{-0.15cm}$$\mb{S}\lra L_{\br}^*\hra D_{\br}^*\lra G_{\br}$$ where the first map is the inclusion of the coordinate corresponding to $\phi_{2n+1}$ and the last map is the norm.

For $h' = h|_{S^1} :S^1\lra G'_{\br}$, there is a map $\ti{h'}: S^1\lra (T'_L)_{\br}\sub {G'_{\br}}$ such that $h' = \rho'\circ \ti{h'}$. Now $h$ is conjugate to the map $h_0$ and factors through the torus $T\sub G$, the image of $L^*\sub D^*$. So $h$ defines a set of special points of the Mumford curve. An abelian variety in the corresponding isogeny class is of dimension $2^{2n}$ and has complex multiplication. At any place $v$ of good ordinary reduction with no simple components of multiplicity larger than one, we have $\End^0(X) = \End^0(X_v)$. So to study the splitting of $X_v$, it suffices to study that of $X$. 

\begin{Prop} With the setup as above, $X$ is isogenous to a product $X^{(1)}\times X'$ where $X^{(1)}$ is an ordinary elliptic curve with CM by the field $E$.\end{Prop}
 
\begin{prf} Write $2^{2n}-1 = m(2n+1)$. We define the map $\tau: T_L'\lra T_L'$ by $x\mapsto N_{L/E}(x)x^{-2}$. This yields the map \vspace{-0.15cm}$$\pi'=(\tau^{(m)},N_{L/E}):T_L'\lra T_L'^{(m)}\times T_E'.$$ The natural action of $T_L'^{(m)}\times T_E'$ on $L^{(m)}\oplus E\cong \bq^{2^{2n+1}}$ makes $\pi'$ into a $2^{2n+1}$-dimensional representation of $T_L'$. The weights of the representation induced by $N_{L/E}$ are $\pm \left(([\phi_1]+\cdots+[\phi_{2n+1}])\right)$ and those of the representation induced by $\tau^{(m)}$ are $\sum\limits_{i\in \mc{I}}[\phi_i]-\sum\limits_{j\in \mc{J}}[\phi_j]$, where $\mc{I}, \mc{J}$ run through the pairs of non-empty sets $\mc{I}, \mc{J}$ such that $\mc{I}\cup \mc{J}=\{1,\cdots,2n+1\}$, $\mc{I}\cap \mc{J}=\emptyset$. So the weights of $\pi'$ are $\sum\limits_{i=1}^{2n+1} \pm[\phi_i]$. These coincide with the weights of the representation $\rho' = N\circ \ti{\rho'}:T_L'\rightarrow G'$ and hence, $\pi'=\rho'$. Thus, $X_v$ is isogenous to a product $X^{(1)}\times X'$ where $X^{(1)}$ is an ordinary elliptic curve with CM by $E$ and $X'$ is an abelian variety (not necessarily simple) with an action of $L$. \end{prf}
\vspace{-0.2cm}
\noindent\underline{Remark} The $2^{2n}-1$-dimensional complement $X'$ is not simple when $n\geq 2$. The Galois group $\Gal(\W_{X_v}/\bq)$ is $\{\pm 1\}\times H$ for some transitive $H\sub S_{2n+1}$ and hence, the simple decomposition of $X'$ is determined by $H$.

The following proposition now implies the existence of generic points on the Mumford curve such that the corresponding abelian variety has a reduction at some prime that has an ordinary elliptic curve as a simple component.

\begin{Prop} With the setup as above, let $x_F$ be a special point on the Mumford curve $\mc{M}/F$ and $X$ the corresponding abelian variety with complex multiplication. Let $\mf{p}$ be a prime of $F$ with residue field $k$. There exists a finite extension $F'$ of $F$, a prime $\mf{p'}$ lying over $\mf{p}$ and a generic point $y_{F'}\in \mc{M}(F')$ such that the abelian variety $Y$ corresponding to $y_{F'}$ has reduction $Y_{k'}$ isogenous to $X_{k'}$.\end{Prop}
\vspace{-0.2cm}
The proof is identical to that of ([No01], Proposition 5.1) and hence, is omitted. We conclude this section by briefly describing the ordinary elliptic curves over finite fields that occur as simple components of the reductions $X_v$. 

As shown in [SZ12], the reduction $A_v$ is ordinary if and only if the place $v$ has local degree $1$ in the field extension $K/\bq$. Since $K/\bq$ is Galois, this is equivalent to $p$ splitting completely in $K$. As seen above, the elliptic curve $X^{(1)}$ has CM by the imaginary quadratic field $E$. Since $\End^0(X^{(1)})$ determines $X^{(1)}$ up to isogeny, it suffices to determine the possibilities for the endomorphism algebra $\End^0(X^{(1)})$.

Let $\mf{l}_1,\cdots,\mf{l}_t$ be the finite primes of $K$ that $D$ is ramified at and let $l_i$ be the rational prime lying under $\mf{l}_i$. Now $EK$ has an embedding in $D$ if and only if $EK$ splits $D$. Since $EK$ is a CM field, this is equivalent to $EK$ having local degree $2$ at all the finite primes of $K$ that $D$ is ramified at. But since we assumed that $K/\bq$ is Galois of odd degree, this is equivalent to $E/\bq$ being either ramified or inert at each of the rational primes $l_i$.

In appendix 6.3, we show a way to construct abelian varieties strongly of Mumford type with simple reduction at some finite place (Proposition 6.3).

\vspace{-0.2cm}
\begin{center}
\section{\fontsize{11}{11}\selectfont Density of places of simple reduction} 
\end{center}
\vspace{-0.2cm}
In this section, we study the density of finite places where the different types of reductions occur. Since the abelian varieties weakly of Mumford type have endomorphism rings $\bz$, the following conjecture suggests that they have absolutely simple reduction at a set of places of density one.

\noindent \textbf{Conjecture} ([MP08], [Zyw14]) \textit{Let $X$ be an absolutely simple abelian variety over a number field $F$ such that $F = F_X^{\conn}$. There exists a density one set of places of $F$ such that the reduction $X_v$ is isogenous to the $d$-th power of a simple abelian variety where $d^2$ is the dimension of the division algebra $\End^0_{\qbar}(X)$ over its center.}

While the conjecture remains open in this generality, several special cases have been studied. If the abelian variety is, in fact, \textit{strongly} of Mumford type, the conjecture holds by the following theorem from [Zyw14].

\begin{Thm} $\mathrm{[Zyw14]}$ Let $X$ be an abelian variety over a number field $F$ such that $F_X^{\conn} = F$. If $\End^0(X_{\qbar})$ is of dimension $d^2$ over its center and the Mumford-Tate conjecture holds for $X$, then $X_v$ is isogenous to the $d$-th power of an absolutely simple abelian variety for all $v$ outside a density zero set of places of $F$. \end{Thm}
\vspace{-0.2cm}
From this, one can easily deduce our main result for abelian varieties \textit{strongly} of Mumford type.

\begin{Prop} Let $A$ be an abelian variety \textit{strongly} of Mumford type over a number field $F$ such that $F_A^{\conn} = F$. Then $A$ has absolutely simple reduction at a set of places of density one. \end{Prop}

\begin{prf} Since the root system of $\Lie(\MT(A))$ is of type $A_1$, the Mumford-Tate conjecture holds for $A$ by ([Pin98], Theorem 5.13). Since $\End(A_{\qbar}) = \bz$, this result now follows immediately from Zywina's theorem.\end{prf}
\vspace{-0.1cm}
We note that for dimension $4$ (equivalently, rank $4$), this was proved in ([Ach12], Theorem C). We now study the reductions of abelian varieties weakly of Mumford type. While the Mumford-Tate conjecture remains open for this class for rank greater than $4$ (alternatively dimension greater than $4$), we know from [Pin98, 5.13] that there exists a connected reductive group that looks like the Mumford-Tate group is expected to.
\vspace{-0.2cm}
\subsection{\fontsize{11}{11}\selectfont Reductions of abelian varieties weakly of Mumford type}
\vspace{-0.2cm}
Throughout the rest of this section, unless otherwise specified, let $A$ be an abelian variety of rank $N+1$ over a number field $F$, weakly of Mumford type. Some of the results (Proposition 5.11) generalize to absolutely simple abelian varieties with trivial endomorphism rings and we hope to use the same techniques to study this broader class of abelian varieties. Enlarging $F$ if necessary, we assume all the $G_l$ are connected. By Theorem 5.11 of [Pin98], there exists a connected reductive group $\GA$ over $\bq$ with a faithful absolutely irreducible $\bq$-linear representation $\VA$ such that 

-$(\GA,\VA)\otimes \bq_l\cong (G_l,V_l)$ for all primes $l$ in a set $\mathbb{S}_A$ of density one. 

-$\GA^{\der}$ is $\bq$-simple.

\noindent\underline{Notation.} This section is heavily reliant on adapting the techniques from [Zyw14] and [No09]. To avoid confusion with the notation, we emphasize that while in [Zyw14] the group $\GA$ is the Mumford-Tate group of $A$, we use the symbol for the reductive group over $\bq$ described above. The two groups are the same if we assume the Mumford-Tate conjecture for $A$ (which remains open for our class of abelian varieties, as far as we know).
\vspace{0.05cm}

\noindent \underline{The set $\mc{S}_A$.} We define the set $\mc{S}_A$ of the places $v$ of $F$ such that:

\noindent $1$. $A$ has good ordinary reduction at $v$.\\
$2$. $v$ has local degree $1$ over $\bq$.\\
$3$. $\Phi_{A_v}$ is a free abelian group with $\rank \Phi_{A_v}=\mr{rank}(A)$.

The first condition holds for a set of places of density one by ([Pin98], Theorem 7.13). The second holds for a set of places of density one by the Chebotarev density theorem. That the third condition holds for a density one set is an implication of the following result from [LP97]. For every place $v$, fix a semisimple element $t_v$ in the conjugacy class of $\rho_{_l}(\Fr_v)$ and let $\bT_{v,l}$ be its Zariski closure in $G_l$.

\begin{Thm} $\mr{([LP97],\; Theorem\;1.2)}$ For any $l$, there is a Zariski closed subvariety $Y\subseteq G_l$ such that $\bT_{v,l}$ is a maximal torus of $G_l$ whenever $\rho_{_l}(\Fr_v)\in G_l\setminus Y$.\end{Thm}
\vspace{-0.2cm}
Now $\bT_{v,l}$ is a maximal torus if and only if $\Phi_{A_v}$ is a free abelian group of the same rank as the rank of $A$. Thus, the third condition holds for a set of places of density $1$ and hence, $\mc{S}_A$ has density $1$. So, to prove Theorem 1.1, it suffices to show that for some density zero subset $\mc{V}$ of $\mc{S}_A$, $A$ has absolutely simple reduction at every $v\in \mc{S}_A\setminus \mc{V}$.

We note that for any $v\in \mc{S}_A$ and prime $l\in\mb{S}_A$ other than $\charac(v)$, the element $\rho_{_l}(\Fr_v)$ is semisimple and hence, lies in some maximal torus of $G_l(\bq_l) = \GA(\bq_l)$. So for any maximal torus $\bT$ of $\GA$, $\bT(\qbar_l)$ contains an element $t$ conjugate to $\rho_{_l}(\Fr_v)$ in $G_l(\qbar_l) = \GA(\qbar_l)$. The following proposition shows that good ordinary reduction implies $\Phi_{A_v}$ being of full rank if $N$ is a prime and $A$ is an abelian variety weakly of Mumford type of rank $N+1$.

\begin{Prop} Suppose the integer $N = \mr{rank}(A)-1$ is a prime. Let $v$ be a place of good ordinary reduction. Then the map $\ev :X^*(\bT)\lra \Phi_{A_v}$ is an isomorphism.\end{Prop}
\begin{prf} Since the map is surjective for all places of good reduction, it suffices to show  injectivity. But since $N$ is a prime, appendix 6.2 shows that $\ker(\ev)$ is one of the following:

\noindent (a) $0$\\
(b) $\bz^{N}$\\
(c)$\{ (x,\cdots,x):x\in\bz\}$\\
(d)$\{(x_{_1},x_{_2},\cdots,x_{_N}:\sum\limits_{i=1}^N x_i=0\}$

So it suffices to verify that cases (b), (c), (d) do not occur for places of ordinary reduction. In case (b), $A_v$ is a power of the supersingular elliptic curve. In case (c), $\ev(1,\cdots,1)\in \ker(\ev)$ $\ev(1,\cdots,1)\sqrt{q}$ corresponds to a supersingular elliptic curve. In case (d), $\rho_{_l}(\Fr_v)$ has precisely $N+1$ distinct eigenvalues. Hence, $\ev$ is injective for places of ordinary reduction, which completes the proof.\end{prf}
\vspace{-0.15cm}
Let $\brho_{_l}:\absf\lra \Aut_{_{\bz/l\bz}}(A[l])$ be the representation $\si\mapsto (P\mapsto P^{\si})\;\forall\;P\in A[l]$. Then $\brho_{_l}(\absf)$ injects into $\G_l(\bF _l)$. The following proposition shows that it is not much smaller than the full $\G_l(\bF _l)$.

\begin{Prop} $\mathrm{[Ser00]}$ Suppose $F_{X}^{\conn} = F$ and $X$ is absolutely simple.\\
$1$. For $l$ sufficiently large, $\G _l$ is a reductive group over $\bz_l$.\\
$2$. There is a constant $C$ such that $[\G _l(\mb{F}_l):\brho_{_l}(\absf)]\leq C$.\\
$3$. For $l$ sufficiently large, the group $\brho_{_l}(\absf)$ contains the commutator subgroup $\G_l(\mb{F}_l)$.\end{Prop}
\vspace{-0.2cm}
A sketch of the proof may be found in [Zyw14]. More details may be found in [Ser00] and [Win02]. Combining all the $l$-adic representations attached to $A$, we define the representation \vspace{-0.1cm}
$$\rho_{_A}:\absf\lra \prod\limits_l \Aut_{\bq_l}(V_l(A)),\;\;\;\;\si\mapsto (\rho_{_l}(\si))_l.$$
After enlarging the field $F$ if necessary, the representations $\rho_{_l}$ are independent, by the following result of Serre.

\begin{Prop} $\mr{([Ser00], 138)}$ There exists a finite extension $F'/F$ such that $\rho_{_A}(\Gal _{F'}) = \prod\limits_l \rho_{_l}(\Gal _{F'})$.\end{Prop}
\vspace{-0.2cm}
This field extension will be important in measuring the density of the places of $F$ whose Frobenius lands in a fixed conjugacy class of $\brho_{_{l}}(\absf)$.  

\begin{Prop} $\mr{([Zyw14], \mr{\;Proposition\;}\; 2.12)}$ Let $F'/F$ be an extension as in Serre's lemma. Let $\Lamb$ be a finite set of rational primes. For each $l\in \Lamb$, fix a subset $U_l$ of $\brho_{_l}(\absf)$ stable under conjugation. The set of places $v$ of $F$ such that $\brho_{_l}(\Fr_v)\subseteq U_l\;\forall\; l\in \Lamb$ has density \vspace{-0.15cm}$$ \sum\limits_{C}\frac{|C|}{\Gal(F'/F)}\prod\limits_{l\in\Lamb}\frac{|\brho_{_l}(\Gamma_C)\cap U_l|}{|\brho_{_l}(\Gamma_C)|}$$ where $C$ varies over the conjugacy classes of $\mr{Gal}(F'/F)$ and $\Gamma_C:=\{\sigma\in \Gal_F:\sigma|_{F'}\in C\}$. \end{Prop}

\underline{The groups $W(G,T)$ and $\Pi(G,T)$:} Let $G$ be a reductive group over a perfect field $K$. Its maximal tori are all conjugate to each other over the algebraic closure. Choose a maximal torus $T$ of $G$. We have a homomorphism \vspace{-0.1cm} $$\phi_{_T}:\absk\lra \Aut(X^*(T)),\;\;\;\;\; \si\mapsto\; (\chi \mapsto \si(\chi))$$ for $\si\in \absk$, $\chi\in X^*(T)$. Note that $T$ is split if and only if $\phi_{_T}(\absk) = \{1\}$.

We use the symbol $W(G,T)$ to mean the Weyl group $N_{G}(T)(\ov{K})/T(\ov{K})$ where $N_{G}(T)$ is the normalizer of $T$ in $G$. Now, $\absk$ acts on $W(G,T)$ and $\si(w) = \phi_{_T}(\si)\circ w\circ \phi_{_T}(\si)^{-1}$ for $\si\in \absk$, $w\in W(G,T)$. We define $\Pi(G,T)$ to be the subgroup of $\Aut(X^*(T))$ generated by $W(G,T)$ and $\phi_{_T}(\absk)$.

In particular, if $A$ is an abelian variety weakly of Mumford type of rank $N+1$ and $\GA$ the reductive group over $\bq$ constructed in [Pin98, Theorem 5.13], we have $W(\GA,\bT)\cong \{\pm 1\}^{N}$. Let $\Psi$ be a root system associated to $\GA$. Since $\GA^{\der}$ is semisimple, the group $\Pi(\GA)$ is isomorphic to a subgroup of $\Aut(\Psi)$ which contains the Weyl group $W(\Psi) = W(\GA)\cong \{\pm 1\}^{N}$. The group $\Aut(\Psi)$ is isomorphic to the semidirect product $\{\pm 1\}^{N}\rtimes S_{_N}$ where $S_{_N}$ acts on $\{\pm 1\}^{N}$ by permuting  the coordinates. So $\Pi(\GA)$ is isomorphic to $\{\pm 1\}^{N}\rtimes \phi_{_{\bT}}(\absq)$ where $\phi_{_{\bT}}(\absq)$ is some transitive subgroup of $S_{_N}$.

Let $v$ be a place of good reduction and $l$ any prime other than $p:=\charac(v)$. Since $A$ has good reduction at $v$, there exists a polynomial $P_{A_v}(X)\in \bz[X]$ such that $P_{A_v}(X) = \det(XI-\rho_{_l}(\Fr_v))$ for all $l\neq p$. For each character $\chi\in X^*(\bT)$, define $\VA(\chi) = \{v\in \VA\otimes_{\bq} \qbar: tv = \chi(t)v\;\forall\;t\in \bT(\qbar)\}$. Define the set of weights $\Om = \{\chi\in X^*(\bT): \VA(\chi)\neq 0\}$. We have the decomposition $\VA\otimes_{\bq} \qbar = \bigoplus\limits_{\chi\in \Om} \VA(\chi)$ and hence, \vspace{-0.15cm}$$\det(XI-t) = \prod\limits_{\chi\in \Om} (X-\chi(t)).$$ The set $\Om$ is stable under the actions of $W(\GA,\bT)$ and $\absq$ on $X^*(\bT)$ which gives a natural action of $\Pi(\GA,\bT)$ on $\Om$. The representation $\VA\otimes_{\bq}\qbar$ of $\GAQ$ is irreducible since $\End^0(A_{\qbar}) = \bq$ and $W(\GA,\bT)$ acts transitively on $\Om$. The following result will be crucial for the density one result.

\begin{Thm} Let $A$ be an abelian variety weakly of Mumford type of rank $N+1$ over a number field $F$ such that $F_A^{\conn} = F$. For any number field $L$ splitting $\GA$, we have $\Gal(L(\W_{A_v})/L) = W(\GA,\bT)$ away from a set of places of density zero. \end{Thm}

\vspace{-0.2cm}
The proof will cover the rest of this section. We will draw  heavily from the techniques from [Zyw14] and [No09], although we will not assume the Mumford-Tate conjecture for our class of abelian varieties.
\vspace{-0.4cm}
\subsection{\fontsize{11}{11}\selectfont Frobenius conjugacy classes}
\vspace{-0.2cm}
\underline{The variety $\Conj'(\GA)$.} Let $G\sub \GL_{2g,\bq}$ be a connected reductive algebraic group over $\bq$ such that its center is one-dimensional and $G^{\der}$ is $\bq$-simple. For simplicity, we assume that $\Lie(G^{\der})$ is $\otimes$-isotypic and the simple factors are of type $A$, $B$, $C$ or $D^{\mb{H}}$. Let $G_1,\cdots,G_d$ be the simple factors of $G_{\qbar}^{\der}$. We have a central isogeny\vspace{-0.2cm} $$G_1\times \cdots \times G_d\lra G_{\qbar}^{\der},\;\;\;\; (g_1,g_2,\cdots,g_d)\mapsto g_1g_2\cdots g_d .$$ Let $\mc{A}$ be the group $f$ of automorphisms of $G_{\qbar}$ such that\\
- $f$ preserves each $G_i$.\\
- the restriction $f|_{G_i}:G_i\lra G_i$ agrees with conjugation by some element of $G_i$.\\
- $f$ acts as the identity on the one-dimensional center of $G_{\qbar}$

Let $R$ be the affine coordinate ring of $G$, on which $G$ acts by composition with inner automorphisms. The group $\mc{A}$ acts on $R$ by composition. We define $R^{\mc{A}}$ to be the $\bq$-subalgebra consisting of elements of elements of $R$ fixed by $\mc{A}$. Now define $\Conj'(G):= \Spec(R^{\mc{A}})$. This is a universal categorical quotient. Let $\Cl'_G: G\lra \Conj'(G)$ be the morphism arising from the inclusion $R^{\mc{A}}\hra R$ of $\bq$-algebras. For a fixed maximal torus $T$ of $G$, let $\mc{A}(T)$ be the subgroup of $\mc{A}$ that preserves $T_{\qbar}$. Define $\Gam_G: = \{f|_{T_{\qbar}}:f\in \mc{A}(T)\}$. This is a finite subgroup of $\Aut(T_{\qbar})$ stable under the $\absq$-action. Furthermore, $\Conj'(G)$ is the quotient of $T_{\qbar}$ by $\Gam_A$.

For a faithful $d$-dimensional $\bq$-linear representation $V$ of $G$, we have a morphism \vspace{-0.1cm}$$P_{_V}:G\hra\GL(V)\lra \mb{A}_{\bq}^{d}$$ sending an element to the coefficients of its characteristic polynomial. This morphism factors through a map $\ov{P}_{_V}:\Conj'(G)\lra \mb{A}_{\bq}^{d}$ which is a finite morphism of schemes.

Now let $A$ be an abelian variety over a number field $F$, weakly of Mumford type of rank $N+1$. Enlarging $F$ if necessary, we may assume $F_A^{\conn} = F$. Since $\End(A_{\qbar}) = \bz$, there exists a connected reductive group $\GA$ over $\bq$ and a faithful absolutely irreducible $\bq$-linear representation $\VA$ such that $\GA^{\der}$ is $\bq$-simple and $(\GA,\VA)\otimes \bq_l \cong (G_l,V_l)$ for a set of rational primes of density one. The group $\GA$ fulfills all of the conditions imposed on $G$. We define the variety $\Conj'(\GA)$ as described above. We also choose and fix a maximal torus $\bT$ of $\GA$. Note that in this case, $\Gam_A = W(\GA,\bT)\cong \{\pm 1\}^{N}$.\vspace{-0.2cm}

\subsection{\fontsize{11}{11}\selectfont Weakly neat elements}

\begin{Def}\normalfont Let $B$ be an abelian variety over a number field. Following $\mr{[No09]}$, we say a place $v$ of good reduction is \textit{weakly neat} if none of the quotients of eigenvalues $\pi_{_1}\pi_{_2}^{-1}$ of the Frobenius are non-trivial roots of unity. The place is \textit{neat} if the group $\Phi_{B_v}$ is torsion-free.\end{Def}
\vspace{-0.2cm}

We will need the following straightforward generalization. Let $G$ be a reductive group over a field $K$ of characteristic zero and $V$ any faithful $K$-linear representation of $G$.

\begin{Def}\normalfont A semisimple element $g\in G(\Kbar)$ is \textit{weakly neat} if none of the quotients of the eigenvalues of $g$ are non-trivial roots of unity. It is \textit{neat} if the group generated by the eigenvalues is torsion-free. \end{Def}
\vspace{-0.2cm}
The preceding definition is independent of the choice of the faithful representation $V$. Clearly, a place $v$ of a number field being weakly neat with respect to an abelian variety over the number field is equivalent to $\rho_{_l}(\Fr_v)$ being weakly neat for any $l$ other than $p = \charac(v)$.

\underline{Remark.} Let $B$ be an abelian variety $B$ over a finite field $k$. Suppose $\W _B$ has two eigenvalues $\pi_{_1}$, $\pi_{_2}$ such that $(\pi_{_1}\pi_{_2}^{-1})^N=1$. Then $\pi_{_1}^N= \pi_{_2}^N$. Let $\ti{v}$ be a place lying over $v$ in a finite extension over $F$ such that the residue degree is $N$. Then $\W _{A_{\ti{v}}} = \W_{A_{v}}^N$ has the eigenvalue $\pi_{_1}^N$ with multiplicity at least $2$. Thus, for any abelian variety over a number field, all places of good ordinary reduction with $\End(A_v)$ commutative are weakly neat.  In particular, for an abelian variety weakly of Mumford type, this condition holds away from a density zero set by [Pin98] and [LP97].
\vspace{-0.2cm}
\subsection{\fontsize{11}{11}\selectfont A variant of a theorem of Noot}
\vspace{-0.1cm}
We recall (Theorem 2.4) that we have the inclusion $G_l^0\sub \MT(A)\times_{\bq} \bq_l$ of algebraic groups and in particular, we have the inclusion $G_l^0(\bq_l)\sub \MT(A)(\bq_l)$ of $\bq_l$-points for every prime $l$. Replacing the field $F$ by the extension $F_{A}^{\conn}$, we have the inclusion $G_l(\bq_l)\sub \MT(A)(\bq_l)$ for every $l\in\mb{S}_A$. The following compatibility theorem is pivotal in the proof of Theorem 5.1 (the main theorem in [Zyw14]).

\begin{Thm} $\mr{([No09],\; Theorem\; 1.8)}$ Let $A$ be an abelian variety over a number field $K = K_A^{\conn}$. Let $v$ be a weakly neat place of $K$. Then there exists an element $\Cl_v'\in \Conj'(\MT(A))(\bq)$ such that $\Cl_v' = \mr{Cl}'_{\MT(A)}(\rho_{_l}(\Fr_v))$ for all primes $l$ such that $v\nmid l$.\end{Thm}
\vspace{-0.2cm}
Our approach will diverge here from that of [No09] and [Zyw14] since the Mumford-Tate conjecture remains open for abelian varieties weakly of Mumford type of rank greater than $4$ (as far as we know). Instead, we will need a variant of this compatibility theorem that replaces $\MT(A)$ by the group $\GA$ and the set of all primes other than $p=\charac(v)$ by the density one subset $\mb{S}_A\setminus \{p\}$.

For any integer $n\geq 1$, let $\tau_n:\GA\lra \GA$ be the $n$-th power map. Since $\Conj'(\GA)$ is a universal categorical quotient, there exists a unique map $\ov{\tau}_n:\Conj '(\GA)\lra \Conj '(\GA)$ making the following diagram commute.
\vspace{-0.2cm}
\begin{center}
$\begin{CD}
\GA @>{\tau_n}>> \GA\\
@VV{\Cl'_{\GA}}V @VV{\Cl'_{\GA}}V\\
\Conj '(\GA) @>{\ov{\tau}_n}>> \Conj '(\GA) 
\end{CD}$
\end{center}
We note that this is a finite map of $\bq$-varieties, which will be relevant later. For the next two lemmas, we drop the assumption that $A$ is weakly of Mumford type and only assume the following:

\noindent 1. $\End(A_{\qbar}) = \bz$. So the triple $(\GA,\VA,\mb{S}_A)$ is well-defined and $\GA^{\der}$ is $\bq$-simple.\\
2. The root systems of $\Lie(\GAQ^{\der})$ are of type $A$, $B$, $C$ or $D^{\mb{H}}$.

Since the root systems of $\Lie(\GAQ^{\der})$ are of type $A,B,C$ or $D$ for any abelian variety with endomorphism ring $\bz$ ([Pin98]), the second condition is equivalent to assuming that the root systems are not of type $D^{\br}$. We will need the following lemma.

\begin{Lem} Let $\al$, $\be$ be two weakly neat elements of $\GA(\qbar)$ with the same characteristic polynomial in $\VA$. If $\ov{\tau}_n(\al) = \ov{\tau}_n(\be)$, then $\ClGA(\al) = \ClGA(\be)$. \end{Lem}
	
\begin{prf} We write $g=\dim A$. The representation $\VA$ is $2g$-dimensional and the symmetric group $S_{2g}$ acts on $\bg_{m,\qbar}^{2g}$ by permutation. Define the set
\vspace{-0.1cm} $$\bT_w: = 	\{t\in \bT: t\text{ is weakly neat and } \Cent_{S_{2g}}(t) = \Cent_{S_{2g}}(\bT) \}.$$
Let $s,t\in \bT_w$ be elements such that $\tau_n(s) = \tau_n(t)$ and $s,t$ have the same characteristic polynomial. Then there exists some $\tau\in S_g$ such that $\tau(s) = t$. Thus, $\tau(s^n) = t^n$ and hence $\tau\in \Cent_{S_{2g}}(\bT)$. So $\tau$ is the identity on $\bT$ and hence, $t_1 = t_2$. Since $\al$ is a semisimple element in the reductive group $\GA$, it lies in a maximal torus. Write $\al = (\al_1,\cdots,\al_{2g})\in \bg_m^{2g}$. Define the torus \vspace{-0.1cm}$$\bT^{\al} = \{(t_1,\cdots, t_{2g})\in \bg_m^{2g}:t_i = t_j\text{ if } \al_i=\al_j \}.$$ Then, in particular, $\al\in \bT_w^{\al}$. Since $\ClGA(\al^n) = \ClGA(\be^n)$, it follows that $\al^n$ and $\be^n$ are conjugate in $\GA(\qbar)$. So there exists a weakly neat semisimple element $\ga\in \GA(\qbar)$ such that $\al^n = \ga^n$ and $\ClGA(\be) = \ClGA(\ga)$ and hence, $\al$ and $\ga$ have the same characteristic polynomial. Since $\al^n = \ga^n$, we have $\ga\in T_w^{\al}$. Thus, $\ClGA(\al) = \ClGA(\be)$.\end{prf}

In particular, the map\vspace{-0.15cm}$$\ov{\tau}_n\times \ov{P}_{_{\VA}}: \Conj'(\GA)\lra \Conj'(\GA)\times \mb{A}_{\bq}^{2g}$$ is injective on the subset of $\Conj'(\GA)(\qbar)$ that is the image under $\ClGA$ of the subset of weakly neat elements in $\GA(\qbar)$. We now use this lemma to show that the statement of the theorem is true if it is ``potentially'' true, in the sense of passing to a finite extension.

\begin{Prop} Let $A$ be an abelian variety over a number field $F$ such that $F_A^{\conn} = F$ and $\End(A_{\qbar}) = \bz$. Let $v$ be a weakly neat place of $F$. Suppose there exists a finite extension $F'/F$ with a place $\ti{v}$ lying over $v$ such that for some conjugacy class $\Cl'_{\ti{v}}\in \Conj'(\GA)(\bq)$, we have $\Cl'_{\GA}(\rho_{_l}(\Fr_{\ti{v}})) = \Cl'_{\ti{v}}$ for all $l\in \mb{S}_A$. Then there exists a conjugacy class $\Cl'_{v}\in \Conj'(\GA)(\bq)$ such that $\Cl'_{\GA}(\rho_{_l}(\Fr_v)) = \Cl'_{v}$ for all $l\in \mb{S}_A$.\end{Prop}

\begin{prf} We choose a prime $l$ in $\mb{S}_A$. Pick an element $x:= \Cl'_{\GA}(\Fr_v)\in \Conj'(\GA)(\bq)\subseteq \Conj'(\GA)(\qbar_l)$. Let $n$ be the residue degree of $\ti{v}$ over $v$. Let $\ti{x}$ be a pre-image of $\Cl'_{\GA}(\Fr_v)$ under the map $\GA(\qbar_l)\xra{\ClGA} \Conj'(\GA)(\qbar_l)$ of $\qbar_l$-valued points induced by $\ov{\tau}_n$. We first show that $\ti{x}$ is a $\qbar$-valued point of $\Conj'{(\GA)}(\qbar_l)$. 

Let $R^{\mc{A}}$ be the coordinate ring of $\Conj'{(\GA)}$ and let $\ov{\tau}_n^{\sharp}:R^{\mc{A}}\lra R^{\mc{A}}$ be the homomorphism induced on the coordinate rings. Then $\tau_n\circ \ti{x}=x$. So we have the commutative diagram:
\vspace{-0.2cm}
\begin{center}
$\begin{CD}
R^{\mc{A}} @>{\ov{\tau}_n^{\sharp}}>> R^{\mc{A}}\\
@VV{x}V @VV{\ti{x}}V\\
\bq @>{i}>> \qbar_l\\ 
\end{CD}$
\end{center}

Since $\ov{\tau}_n^{\sharp}$ is a finite homomorphism of algebras, $\ti{x}(R^{\mc{A}})$ is a finite algebra over $x(R^{\mc{A}})= \bq$. So the fraction field $K:=\mr{Frac}(\ti{x}(R^{\mc{A}}))$ is a number field. The diagram above may be replaced by the following diagram:
\vspace{-0.4cm}
\begin{center}
$\begin{CD}
R^{\mc{A}} @>{\ov{\tau}_n^{\sharp}}>> R^{\mc{A}}\\
@VV{x}V @VV{\ti{x}}V\\
\bq @>{i}>> K\\ 
\end{CD}$
\end{center}

Since $\ti{x}$ is a $K$-valued point of $\Conj'{(\GA)}$, this proves the claim that $\ti{x}$ is a $\qbar$-valued point of $\Conj'(\GA)$.

Now $\Cl'_{\GA}(\rho_{_l}(\Fr_v))$ is a $\bq_l$-valued point of $\Conj'{(\GA)}$ for every  $l(\neq p)$ in $\mb{S}_A$. Let $l$, $l'$ be any two distinct primes in $\mb{S}_A$ other than $p$. Then we have a weakly neat element $\rho_{_{l'}}(\Fr_v)$ with the same characteristic polynomial as that of $\rho_{_l}(\Fr_v)$. By the hypothesis, we have a unique element $\Cl_{\ti{v}}'\in \Conj'{(\GA)}(\bq_l)$ such that \vspace{-0.1cm}$$\ov{\tau}_n(\Cl'_{\GA}(\rho_{_l}(\Fr_v))) = \ov{\tau}_n(\Cl'_{\GA}(\rho_{_{l'}}(\Fr_v))) = \Cl_{\ti{v}}'.$$ By the preceding discussion, $\Cl'_{\GA}(\rho_{_{l}}(\Fr_v))\in \Conj'(\GA)(\qbar)$. So from lemma 5.10, it follows that $\Cl'_{\GA}(\rho_{_l}(\Fr_v)) = \Cl'_{\GA}(\rho_{_{l'}}(\Fr_v))$ for all $l,l'\in \mb{S}_A$ other than $p$. We denote this common element by $\Cl_v'$. 

It remains to show that $\Cl_v'$ is $\bq$-valued rather than simply $\qbar$-valued. Now, $\Cl_v'\in \Conj'(\GA)(\bq_l)$ for all $l\in \mb{S}_A$ other than $p$. Suppose, by way of contradiction, that $\Cl_v'\notin \Conj'(\GA)(\bq)$ and let $K(\neq \bq)$ be the smallest number field such that $\Cl_v'\in \Conj'(\GA)(K)$. Let $\al\in \mc{O}_K$ be a primitive element for $K/\bq$ and let $f\in \bz[X]$ be the irreducible monic polynomial of $\al$ over $\bq$ with $\deg(f)\geq 2$. By the Chebotarev density theorem and Jordan's lemma, the set of rational primes \vspace{-0.1cm}$$\mc{P}:= \{l: \text{ The image of } f(X) \text{ in } \mb{F}_l[X] \text{ is separable and has no zero} \}$$ has positive density. For any $l\in \mc{P}$, Hensel's lemma implies that $K$ does not have an embedding in $\bq_l$. Since $\mb{S}_A$ is of Dirichlet density one and is disjoint from $\mc{P}$, we have a contradiction. Thus, $\Cl_v'\in \Conj'(\GA)(\bq)$. \end{prf}
\vspace{-0.05cm}
\begin{Prop} Let $A$ be an abelian variety weakly of Mumford type of rank $N+1$ over a number field $F = F_A^{\conn}$. Let $(\GA,\VA,\mb{S}_A)$ be a triple as in {\normalfont ([Pin98], Theorem 5.1)}. Let $v$ be a weakly neat place of $F$. Then there exists a conjugacy class $\Cl'_{v}\in \Conj'(\GA)(\bq)$ such that $\ClGA(\rho_{_l}(\Fr_v)) = \Cl'_v$ for all $l\in \mb{S}_A$ other than $p = \charac(v)$.\end{Prop}

\begin{prf} Since $\End(A_{\qbar}) = \bz$, the $\GAQ$-representation $\VA\otimes_{\bq} \qbar$ is irreducible. The group of homotheties of $\VA\otimes_{\bq} \qbar$ commutes with the action of $\Lie(\GAQ^{\der})$ which is isomorphic to $\mf{sl}_{2,\qbar}^{N}$. By the preceding lemma, it suffices to show this result over a finite extension. Enlarging the field if necessary, we may assume the place $v$ is neat. We write $g = 2^{N-1}$ for brevity. 

Define the map ${P}_{_{\VA}}: \GA\lra \mb{A}_{\bq}^{2g}$ sending an element to the coefficients of its characteristic polynomial with respect to $\VA$. It factors through a finite morphism $\ovPV: \Conj'(\GA)\lra \mb{A}_{\bq}^{2g}$. Let $\ti{S}_{w.n.}\sub \GA(\qbar)$ be the subset of weakly neat elements and let $S_{w.n.}$ be its image in $\Conj'(\GA)(\qbar)$ under the map $\ClGA$.

Now $\Lie(\GAQ^{\der})\cong \mf{sl}_{2,\qbar}^{N}$ and in particular, $\Lie(\GAQ^{\der})$ is $\otimes$-isotypic of type $A$. So the characteristic polynomial determines the conjugacy class of a semisimple element. In particular, the map $\ov{P}_{_{\VA}}: \Conj'(\GA)\lra \mb{A}_{\bq}^{2g}$ is injective on the subset $S_{w.n.}$. Hence, there is a unique element $\Cl'_v\in \Conj'(\GA)(\bq)$ such that $\ov{P}_{_{\VA}}(\Cl'_v) = P_{_{\VA}}(\Fr_v)$. Thus, $\Cl'_v = \Cl'_{\GA}(\rho_{_l}(\Fr_v))$ for all $l\in \mb{S}_A$. Now, $\Cl_v'\in \Conj'(\GA)(\bq_l)$ for all $l\in\mb{S}\setminus\{p\}$. Since $\mb{S}_A$ is of density one, we have $\Cl_v'\in \Conj'(\GA)(\bq)$ by the Chebotarev density theorem.\end{prf}

\vspace{-0.6cm}

\subsection{\fontsize{11}{11}\selectfont Galois representations}
\vspace{-0.3cm}

We return to the setting where $A$ is an abelian variety over a number field $F = F_A^{\conn}$ weakly of Mumford type of rank $N+1$. For a place $v$ of good reduction, we have the surjective homomorphism $\ev :X^*(\bT)\lra\Phi_{A_v},\;\; \chi\mapsto \chi(t_v)$ that maps $\Omega$ onto $\W_{A_v}$. For places in $\mc{S}_A$, it is an isomorphism.

\begin{Lem} Let $v\in \mc{S}_A$ and fix an element $t\in \bT(\qbar)$ such that $\det(XI-t) = P_{A_v}(X)$. Then the map \vspace{-0.15cm}$$\ev :X^*(\bT)\lra \Phi_{A_v},\;\;\;\;\chi\mapsto \chi(t) $$ is an isomorphism mapping $\Omega$ onto $\W_{A_v}$.\end{Lem}

\begin{prf}  It is clear that the map is a homomorphism $X^*(\bT)\lra \qbar^*$. The zeros of $P_{A_v}(X)$ are the values $\chi(t)$ with $\chi\in \Om$. So $\W_{A_v} = \ev(\Om)$. The set $\Omega$ generates $X^*(\bT)$ since $\VA$ is a faithful representation of $\GA$. Since $\ev(\Omega)= \W_{A_v}$ and $\W_{A_v}$ generates $\Phi_{A_v}$, it follows that $\ev(X^*(\bT))=\Phi_{A_v}$. Thus, $\ev$ is a surjective homomorphism of free abelian groups. So it suffices to verify that the ranks of the two groups coincide. 

The third condition in the definition of $\mc{S}_A$ ensures that the rank of the group $\Phi_{A_v}$ is $\mr{rank}(A)$. On the other hand, $X^*(\bT)$ is a free abelian group of the same rank as $\GA$. Since $\ev$ is surjective, it follows that it is an isomorphism.\end{prf}
\vspace{-0.1cm}
The upshot is that for any $v\in \mc{S}_A$, if $W(\GA,\bT)\sub \GalW$, then $\absq$ acts transitively on $\W_{A_v}$ and hence, the reduction is absolutely simple. Thus, detecting the places of simple reduction boils down to detecting the places $v$ such that $W(\GA,\bT)\sub \GalW$.

\underline{The elements $\psi_v(\si)$ and $\psi_{v,l}(\si)$}. Fix a place $v\in \mc{S}_A$. By Proposition 5.12, we may choose an element $t_v\in \bT(\qbar)$ such that $\Cl'_{\GA}(t_v) = \Cl_v'\in \Conj '(\GA)(\bq)$. We may choose $t_v$ so that $\det(XI-t_v) = P_{A_v}(X)$. For every $\si\in \absq$, we define $\psi_v(\si)$ to be the unique	element in $\Aut(X^*(\bT))$ such that the following diagram commutes:
\vspace{-0.2cm}
\begin{center}
$\begin{CD}
X^*(\bT) @>{\ev}>> \Phi_{A_v}\\
@VV{\psi_v(\sigma)}V @VV{\sigma}V\\
X^*(\bT) @>{\ev}>> \Phi_{A_v}\\ 
\end{CD}$
\end{center}

On the other hand, for any prime $l\in \mb{S}_A$, we have $W(\GA,\bT)\cong W(G_l,\bT_{\bq_l})$. Furthermore, if $\bTl$ is split and $v\nmid l$, we have a group homomorphism $\psi_{v,l}:\absql\lra W(G_l,\bT_{\bq_l})$ such that $\si(t_{v,l}) = \psi_{v,l}(\si)^{-1}(t_{v,l})$. This map is independent of the choice of the element $t_{v,l}$. 

\underline{Remark.} Let $K$ be a number field splitting $\bT$. Note that if $l$ is split completely in $K$, then $K_w = \bq_l$ for any prime $w$ of $K$ lying over $l$. Hence, $\bT_{\bq_l}$ is split. Since $\mb{S}_A$ has density one,  the set of primes $l\in \mb{S}_A$ such that $\bT_{\bq_l}$ is split has positive density.

We first verify that the elements $\psi_v(\si)$ and $\psi_{v,l}(\si)$ are the same in the case of abelian varieties weakly of Mumford type. In the more general case of an abelian variety with endomorphism ring $\bz$, we can show (with a proof similar to [Zyw14]) that they lie in the same conjugacy class of $\Gamma_A$.

\begin{Lem} Let $l\in \mb{S}_A$ and suppose $\bT_{\bq_l}$ is split. For any $\sigma\in \absql$, we have $\psi_v(\si) = \psi_{v,l}(\si) \in W(\GA,\bT)$.\end{Lem}

\begin{prf} Choose an element $t_{v,l}\in \bT(\qbar_l)$ conjugate to $\rho_{_l}(\Fr_v)$ in $\GA(\qbar_l)$. Then $\ClGA(t_{v,l}) = \ClGA(\rho_{_l}(\Fr_v)) = \Cl_v'$. So there is a unique $\be\in \Gamma_A=  W(\GA,\bT)$ such that $t_{v,l} = \be(t_v)$. We may verify that $\psi_v(\si)= \be^{-1}\psi_{v,l}(\si)\circ\be$. Since $\psi_{v,l}(\si)\in W(\GA,\bT)$ and $W(\GA,\bT)$ is abelian in this particular case, we have $\psi_v(\si) = \psi_{v,l}(\si) \in W(\GA,\bT)$.\end{prf}
\vspace{-0.1cm}
Let $l$ be a prime such that $\G_l$ is a reductive group scheme over $\bz_l$, as is the case when $l$ is sufficiently large. Choose a maximal torus $\mc{T}_l$ in $\G_l$ and let $\bT_l$ be its generic fiber, which is a maximal torus in $G_l$. For any place $v\in \mc{S}_A$ such that $v\nmid l$, we define the conjugacy class \vspace{-0.1cm}$$\mc{I}_{v,l}:=\{t\in \bT(\qbar_l): t \text{ and } \rho_{_l}(\Fr_v) \text{ are conjugate in } G_l(\qbar_l) \}$$ and fix an element $t_{v,l}\in \mc{I}_{v,l}$. The third condition in the definition of the set $\mc{S}_A$ ensures that the group generated by $t_{v,l}$ is Zariski dense in $\bT_{l,\qbar_l}$. So the action of $W(G_l,\bT_l)$ on $\mc{I}_{v,l}$ is simply transitive. For each $\si\in \absql$, there is a unique $\psi_{v,l}(\si)\in W(G_l,\bT_ l)$ such that $\si(t_{v,l}) = \psi_{v,l}(\si)^{-1}(t_{v,l})$. Since $\bT_l$ is split by assumption, the map \vspace{-0.1cm}$$\psi_{v,l}:\absql\lra W(G_l,\bT_{l}),\;\;\;\; \si\mapsto \psi_{v,l}(\si)$$ is a group homomorphism.

\underline{The field $k_{\GA}$.} We define the field $k_{\GA}$ to be the intersection of all algebraic extensions of $\bq$ that split $\GA$. The field $k_{\GA}$ is a number field Galois over $\bq$. We first verify the inclusion $\psi_v(\Gal_{k_{\GA}})\sub W(\GA,\bT)$.

\begin{Lem} The image $\psi_v(\Gal_{k_{\GA}})$ is a subgroup of $W(\GA,\bT)$.\end{Lem}
\begin{prf} It suffices to show that $\psi_v(\Gal_{L})\subseteq W(\GA,\bT)$ for any number field $L$ splitting $\GA$.

Let $L_0$ be a number field that splits $\bT$. Let $\Lambda$ be the co-finite subset of $\mc{S}(L_0|\bq)\cap \mb{S}_A$ of primes $l$ for which $\psi_v$ is unramified at $l$, $v\nmid l$. The torus $\bT_{\bq_l}$ is split for all $l\in \Lambda$. From the preceding lemma, $\psi_v(\Fr_{_l})\in W(\GA,\bT)$. Since the set $\{\Fr_{_l}:l\in \Lambda\}$ generates $\Gal_{L_0}$, it follows that $\psi_v(\Gal_{L_0})\sub W(\GA,\bT)$.

Now let $L$ be a number field splitting $\GA$. Choose a maximal torus $\bT'$ of $\GA$ for which $\bT'\times_{\bq} L$ is split and an element $g\in \GA(\qbar)$ such that $\bT'_{\qbar} = g\bT_{\qbar}g^{-1}$. Define $t_v' := gt_vg^{-1}$. We have $\ClGA(t_v') = \Cl_v'$ and $\det(XI-t_v') = P_{A_v}(X)$. We define the homomorphism $\psi_v':\absq\lra \Aut(X^*(\bT '))$ such that $\si(\al(t_v')) = (\psi_v'\al)(t_v')$ for all $\al\in X^*(\bT)$ and $\si\in \absq$. 

Let $\ga:\bT_{\qbar}\lra \bT'_{\qbar}\;$ be the conjugation by $g$ and let $\ga_*:\Aut(\bT_{\qbar})\lra \Aut(\bT'_{\qbar})$ be the conjugation by $\ga$. Then $\ga_*(W(\GA,\bT)) = W(\GA,\bT')$. From the definition of $\psi_v(\si)$, we verify that $\psi_v(\si) = \ga^{-1}\circ \psi'_v(\si)\circ\ga = \ga_*^{-1}(\psi'_v(\si))$. So $\psi_v(\absl) \subseteq \ga_*^{-1}(W(\GA,\bT')) = W(\GA,\bT)$.\end{prf}

We now show that the reverse inclusion holds away from a density zero set of places of $F$. This will complete the proof of the main result. The following lemma from [Zyw14] will be necessary.

\begin{Lem} $(\mr{[Zyw14],\; \mr{Lemma}\; 5.1})$ Fix a subset $C$ of $W(G_l,\bT_{\bq_l})$ stable under conjugation. There is a subset $U_l$ of $\rho_{_l}(\absf)$ stable under conjugation with the following properties:

\noindent $\mr(a)$. If $v\in\mc{S}_A$ satisfies $v\nmid l$ and $\brho_l(\Fr_v)\subseteq U_l$, then $\psi_{v,l}$ is unramified and $\psi_{v,l}(\Fr_{_l})\subseteq C$.\\
$\mr(b)$. Let $F'$ be a finite extension of $F$ and let $H$ be a coset of $\Gal_{F'}$ in $\absf$. Then we have \vspace{-0.15cm}$$\frac{\big|\ov{\rho}_{_l}(H)\cap U_l\big|}{\big|\ov{\rho}_{_l}(H)\big|} = \frac{\big|C\big|}{W(G_l,\bT_{\bq_l})} + O(1/l)$$ where the implicit constant depends only on $A$ and $F'$.\end{Lem}
\vspace{-0.2cm}
The key ideas in the following proposition are from ([Zyw14], Proposition 6.6). To avoid confusion with the notation, we reiterate that we use the symbol $\GA$ for the group constructed in ([Pin98], Theorem 5.13) rather than the Mumford-Tate group $\MT(A)$.

\begin{Prop} Let $A$ be an abelian variety weakly of Mumford type over a number field $F$ such that $F = F_A^{\conn}$. Let $L$ be a finite extension of $k_{\GA}$. Then $\psi_v(\absl) = W(\GA,\bT)$ for all places $v\in \mc{S}_A$ away from a subset of Dirichlet density $0$.\end{Prop}

\begin{prf} Fix a full rank $\GA$-stable lattice $\mc{L}\sub \VA$. Let $\mc{T}$ be the Zariski closure of $\bT$ in the group scheme $\GL(\mc{L})$. We fix a prime $l_0$ such that $\G_l$ is reductive group scheme over $\bz_l$ for every prime $l\geq l_0$. For every $l>l_0$,  $\mc{T}_{\bz_l}$ is a torus over $\bz_l$ by [LP95]. Hence, $\mc{T}_{\bz_l}$ is a maximal torus of $\G _l$ for any rational prime $l\in \mb{S}_A$ such that $l\geq l_0$. 

Let $F'/F$ be a finite extension fulfilling the condition of Serre's lemma. Let $c$ be an element of $W(\GA,\bT)$. Since $\Gam_A = W(\GA,\bT)\cong \{\pm 1\}^{N}$, the element $c$ is preserved under conjugation by $\Gamma_A$ in this case. Let $\mc{S}(L|\bq)$ be the set of rational primes split completely in $L$. For $l$ in the set $\mc{S}(L|\bq)\;\cap\;\mb{S}_A$ and $\geq l_0$, we may identify $c$ with an element of $W(G_l,\bT_{\bq_l}) = W(\G _l,\mc{T}_{\bq_l})$. Let $\mc{V}$ be the subset of $\mc{S}_A$ of places $v$ such that $\brho_{_l}(\Fr_v)\nsubseteq U_l$, $l\neq p$. Let $\mc{V}_Q$ be the set of places such that $\brho_{_l}(\Fr_v)\nsubseteq U_l$ for all $l <Q$, $l\neq p$. Then $\mc{V} = \bigcap\limits_{Q=1}^{\infty} \mc{V}_Q$. By ([Zyw14], Proposition 2.12), the set $\mc{V}_Q$ has density \vspace{-0.4cm}$$\sum\limits_\mc{C} \frac{|\mc{C}|}{\Gal(F'/F)} \prod\limits_{l\in \mc{S}(L|\bq)\;\cap\;\mb{S}_A \; \cap[l_0, Q]}\frac{\big|\brho_{_l}(\Gamma_{\mc{C}})\cap (\rho_{_l}(\absf)\setminus U_l)\big|}{\big| 	\brho_{_l}(\Gamma_{\mc{C}}) \big|}$$ where $\mc{C}$ varies over the conjugacy classes of $\Gal(F'/F)$ and $\Gamma_\mc{C} = \{\si\in \absf: \si|_{F'}\in \mc{C}\}$. Furthermore, by ([Zyw14], Lemma 5.1 (b)), \vspace{-0.15cm}$$\mr{Density}(\mc{V}_Q)<< \prod\limits_{l\in \mc{S}(L|\bq)\;\cap \;\mb{S}_A\; \cap\;[l_0, Q]}\left(1-\frac{1}{W(G_l,\bT_{\bq_l})}+O(1/l)\right)$$\vspace{-0.2cm} 
$$= \prod\limits_{l\in \mc{S}(L|\bq)\;\cap \;\mb{S}_A\cap\;[l_0, Q]}\left(1-\frac{1}{W(\GA,\bT)}+O(1/l)\right).$$

Now, $\mc{S}(L|\bq)\cap \mb{S}_A$ has density $\frac{1}{[\wti{L}:\bq]}$ where $\wti{L}$ is the Galois closure of $L$ over $\bq$. In particular, the set $\mc{S}(L|\bq)\cap \mb{S}_A$ is infinite and hence, $\mc{V}$ has density zero. Let $v\in \mc{S}_A\setminus \mc{V}$. Choose a prime $l\in \mc{S}(L|\bq)\cap \mb{S}_A$ such that $v\nmid l$ and $\brho_{_l}(\Fr_v)\sub U_l$. Then $\psi_{v,l}$ is unramified at $l$ and $\psi_{v,l}(\Fr_v) = c$ by ([Zyw14], Lemma 5.1(a)). Since $\{c\}$ is stable under conjugation by $\Gamma_A$ and the elements $\psi_{v}(\Fr_{_l}) = \psi_{v,l}(\Fr_{_l})\in W(\GA,\bT)$, it follows that $\psi_{v}(\Fr_{_l}) = c$. Since $c$ was an arbitrary element of $W(\GA,\bT)$, we have $\psi_v(\absl) = W(\GA,\bT)$.\end{prf}

The main result (Theorem 1.1) now follows: 

\begin{prf} Let $L$ be any number field splitting $\GA$. Suppose $v$ lies in the density one subset of $\mc{S}_A$ described in the preceding Proposition. Since $\W_{A_v}$ generates $\Phi_{A_v}$, the restriction $\psi_v|_{\absl}$ factors through an isomorphism $\Gal(L(\W_{A_v})/L)\xra{\cong} W(\GA)$. So $\absq$ acts transitively on $\W_{A_v}$ and hence, $A_v$ is absolutely simple.\end{prf}

We conclude this section by describing the endomorphism algebra $\End^0(A_v)$ for most places.

\begin{Corr} With the setup as in theorem 1.1, $\GalW = \{\pm 1\}^{N}\rtimes H$ for some $H\sub S_{_N}$ and $\End^0(A_v)$ is isomorphic to the fixed field of $H$.\end{Corr}
\begin{prf} Let $v$ be a place of simple ordinary reduction. Then $P_{A_v}(X)$ is irreducible of degree $2^N$. Set $\pi_{_0} = \ev(1,\cdots,1)\sqrt{q}$. Then $\bq(\pi_{_0})$ lies in the fixed field of $H$ and since $[\bq(\pi_{_0}):\bq] = 2^N$, it follows that $\bq(\pi_{_0}) = \bqW^H$.\end{prf}

\begin{center}
\section{\fontsize{11}{11}\selectfont Appendix} 
\end{center}
\vspace{-0.2cm}

We use appendix 6.1 to relate abelian varieties strongly of Mumford type with Mumford curves. In appendices 6.2 and 6.3 we state and prove a few propositions that are original but of an elementary nature.

\subsection{\fontsize{11}{11}\selectfont Generalized Mumford Curves}
\vspace{-0.2cm}
We verify that abelian varieties with Mumford-Tate groups fulfilling the hypothesis of Definition 1.3 arise as generic points on Mumford curves. The proof is identical to that of [No01, Proposition 1.5] for the case $n=1$. We provide a sketch for the reader's convenience. 

\begin{Prop} Let $X$ be a an abelian variety of dimension $2^{2n}$, $V= H^1(X(\bc),\bq)$ the Betti cohomology and $G= \MT(X)$ its Mumford-Tate group. Assume $(G,V)$ is of Mumford type of rank $2n+2$. Then there exist a totally real number field $K$ of degree $N$ and a quaternion $D/K$ such that:\\
$\mr{(i)}$ $D$ is split at exactly one real place\\ 
$\mr{(ii)}$ $\Cor_{K/\bq}(D) = 0$\\
$\mr{(iii)}$ $X$ is isomorphic to a fiber of the family $\mc{X}/M$ over a $\bc$-valued point of $M=M_C$ $\mr{(}$for any sufficiently small open compact subgroup $C\sub G(\mb{A}_f)$ where $\mb{A}_f$ is the ring of finite adeles$\mr{)}$.\end{Prop}

\begin{prf} Let $\tilde{G}\lra G$ be the universal simple cover with $\ti{G}\cong \mb{G}_{m,\bq}\times \ti{G}'$, $\ti{G}'_{\qbar}\cong \SL_{2,\qbar}^{N}$. Now $G^{\der}$ is $\bq$-simple. So $\absq$ acts transitively on the set of factors of the product $\SL_{2,\qbar}^{N}$. Let $H\sub\absq$ be the subgroup fixing the first factor and set $K:=\qbar^H$. Then $K$ is field of degree $N$ and $G'$ is the Weil restriction from $K$ to $\bq$ of a $K$-form of $\SL_2$. By Kneser's theorem ([Kne65], Satz 3), there exists a central division algebra $D$ of degree $4$ over $K$ such that $\tilde{G'} = (D^*)^{\der}$. Since $\ti{G'}$ acts on a $2^{2n+1}$-dimensional $\bq$-vector space, it follows that $\Cor_{K/\bq}(D) = \Mat_{2^{2n+1}}(\bq)$.

Now $V$ carries a symplectic form $\langle\; .\;, .\; \rangle$. Let $h:\mb{S}\lra G_{\br}$ be the morphism defining the Hodge structure on $V\otimes \br$. The symmetric bilinear form $\langle\; .\;,h(i) .\; \rangle$ is positive definite. This implies that if $\mathcal{H}$ is the real form of $G_{\bc}$ corresponding to  the involution $\ad (h(i))$, then  $\mc{H}^{\der}$ is compact.

The projection of $h$ on precisely one factor $\PSL_{2,\bc}$ of $G_{\bc}^{\der}$ is non-trivial. Since $\mc{H}$ is compact, $K$ is totally real and at least $2n$ of the factors of $\ti{G'}_{\br}$ are isomorphic to $\SU_{2,\br}$. Since the corestriction is trivial, it follows that $\ti{G'}_{\br}\cong \SU_{2,\br}^{2n}\times \SL_{2,\br}$. Since $\mc{H}^{\der}$ is compact, it follows that $h$ is conjugate to the map $h_0:\mb{S}\lra G_{\br}$ derived from\vspace{-0.3cm} $$\ti{h}_0: \mb{S}\lra \GSU_{2,\br}^{2n}\times \GL_{2,\br}\cong D^*_{\br},\;\;\; a+b\sqrt{-1}\mapsto \left(1,1,\cdots,1,\begin{pmatrix} a & -b\\
b & a \end{pmatrix}\right).$$ \vspace{-0.2cm} This completes the proof.\end{prf}

\subsection{\fontsize{11}{11}\selectfont Cyclic lattices with torsion-free quotients}
\vspace{-0.3cm}
We use this subsection to clarify some of the things from section 4.4. In particular, we discuss the possibilities for the kernel of the map $\ev :X^*(\bT)\lra \qbar^*$, the image of which is the multiplicative group $\Phi_{A_v}$ generated by the eigenvalues of the Frobenius. We will need the following notion.

\begin{Def} \normalfont A sublattice of $\bz^N$ is \textit{cyclic} if it is invariant under the $N$-cycle $(1\;2\;3\cdots N)$.\end{Def}
\vspace{-0.3cm}
We are primarily concerned with cyclic lattices $\mc{L}\sub \bz^N$ such that the quotient $\bz^N/\mc{L}$ is torsion-free. This will help us deduce the possibilities for $\ker(\ev)$ and consequently, the image of $\ev$. We 
define a map\vspace{-0.15cm}$$\phi:\bz[X]/(X^N-1)\lra \bz^N,\;\;\;\; \sum_{j=0}^{N-1} a_jX^j\mapsto (a_0,a_1,\cdots,a_{_{N-1}})\in \bz^N.$$ For an ideal $J$ of $\bz[X]/(X^N-1)$, the image $\phi(J)$ is a  cyclic sublattice of $\bz^N$. Conversely, for a cyclic lattice $\Lambda$, $\{\sum\limits_{j=0}^{N-1} a_jX^j:(a_1,\cdots,a_{_N})\in \Lambda\}$ is an ideal in $\bz[X]/(X^N-1)$. This gives the bijection\vspace{-0.1cm}$$\{\text{Cyclic sublattices of } \bz^N\} \longleftrightarrow \{\text{Ideals of }\bz[X]/(X^N-1) \}.$$
\vspace{-0.2cm}
Now, $X^N-1= \prod\limits_{d|N} \Phi_d(X)$ where $\Phi_d(X)$ is the $d$-th cyclotomic polynomial  and hence,\vspace{-0.2cm} $$\bz[X]/(X^N-1)\cong \prod\limits_{d|N} \bz[\zeta_d].$$ Since the rings $\bz[\zeta_d]$ are rings with unit elements, an ideal $J$ of $\bz[X]/(X^N-1)$ is a direct product $\prod\limits_{d|N} J_d$ where $J_d$ is an ideal in $\bz[\zeta_d]$. Furthermore, the rings $\bz[\zeta_d]$ are Dedekind domains with finite residue fields and hence, the quotients $\bz[\zeta_d]/J_d$ are finite rings unless $J_d = 0$. So if the quotient is torsion-free, each ideal $J_d$ is either $0$ or the full $\bz[\zeta_d]$. So there are precisely $2^{\tau(N)}$ cyclic sublattices with torsion-free quotient, where $\tau(N)$ is the number of divisors of $N$.

In particular, if $N$ is a prime, the only ideals of $\bz[X]/(X^N-1)\cong \bz\times \bz[\zeta_{_N}]$ with torsion-free quotients are $0$, $\bz\times \bz[\zeta_{_N}]$, $\bz\times 0$ and $0\times \bz[\zeta_{_N}]$. So the only cyclic sublattices $\mc{L}\sub\bz^N$ with a torsion-free quotient $\bz^N/\mc{L}$ are the following:

\noindent (a) $0$\\
(b)$\{ (x,\cdots,x):x\in\bz\}$\\
(c)$\{(x_{_1},x_{_2},\cdots,x_{_N}):\sum\limits_{i=1}^N x_i=0\}$.\\
(d) $\bz^N$

which verifies the claim in Section 4.
\vspace{-0.4cm}

\subsection{\fontsize{11}{11}\selectfont Central division algebras over totally real fields}
\vspace{-0.2cm}
We use this subsection to clarify some of the things in Section 5.6. We will need the following lemma for the proof of Proposition 4.13.

\begin{Lem} Let $d$ be an odd integer. For a finite set of primes ${l_1},\cdots,{l_r}$, there exist infinitely many CM fields $E$ of degree $2d$ over $\bq$ such that:\\
$\mr{(i)}$ $E/\bq$ is cyclic of degree $2d$\\
$\mr{(ii)}$ The primes ${l_1},\cdots,{l_r}$ are inert in $E/\bq$.\end{Lem}

\begin{prf} Choose a prime $\mf{p}$ in $\bq(\zeta_{d})$ such that:\\
- $\mf{p}$ is prime to $2d{l_1}\cdots{l_r}$\\
- $\mf{p}$ has local degree $1$ over $\bq$ \\
- $\mf{p}$ is inert in the degree $2d$ cyclic extension $\bq(\zeta_{4d^2})/\bq(\zeta_{d})$\\
- $\mf{p}$ is split completely in the Kummer extension $\bq\left(\zeta_{2d},\{\sqrt[2d]{l_1l_j}:2\leq j\leq r\}\right)$ but is inert in the cyclic extension $\bq(\zeta_{2d},\sqrt[2d]{l_1})$.

Note that by Kummer theory, $\bq(\zeta_{2d},\sqrt[2d]{l_1})\cap \bq\left(\zeta_{2d},\{\sqrt[2d]{l_1l_j}:2\leq j\leq r\}\right) = \bq\left(\zeta_{2d}\right)$ since\\ $\la l_1\ra \cap \la l_1l_2,\cdots,l_1l_r \ra = \{1\} \sub \bq^*$. Hence, the set of primes $\mf{p}$ of $\bq(\zeta_{d})$ fulfilling the four conditions has positive density by the Chebotarev density theorem. By construction, $\mf{p}$ is inert in each of the cyclic extensions $\bq(\zeta_{d},\sqrt[2d]{l_j})/\bq(\zeta_{d})$ ($1\leq j\leq r$). Let $p$ be the rational prime lying under $\mf{p}$. Then $p\equiv 1\pmod{2d}$ by the second condition and $\gcd(\frac{p-1}{2d}, 2d) =1$ by the third. Define $E$ to be the unique degree $2d$ number field contained in $\bq(\zeta_p)$. Since $\bq(\zeta_p)/\bq$ is cyclic, so is $E/\bq$. Furthermore, since $E$ is of odd index in the CM field $\bq(\zeta_p)$, it is also a CM field. Every rational prime $l_i$ has its inertia degree in $\bq(\zeta_p)$ divisible by $2d$ and since $\gcd(\frac{p-1}{2d}, 2d) =1$, it follows that every $l_i$ is inert in $E/\bq$.\end{prf}
\vspace{-0.1cm}
We now use this lemma to prove Proposition 4.13 from section 4.6.

\noindent\underline{Notations:} For a number field $K$ and a prime $\mfp$, $K_{\mfp^h}$ is the unique degree $h$ unramified extension of the local field $K_{\mfp}$. $\Br(K)$ is the Brauer group of $K$ and for an extension $L/K$, $\Br(L|K)$ is the relative Brauer group, which is the kernel of the restriction map \vspace{-0.1cm}$$\Br(K)\lra \Br(L),\;\;\;\;[D]\mapsto [L\otimes_K D].$$ As before, $\mc{S}(L|K)$ is the set of primes of $K$ split completely in $L$ (equivalently, in the Galois closure of $L$ over $K$). If $L/K$ is Galois, for any place $w$ of $K$, we use the symbol $L^w$ for the completion of $L$ at any place lying over $w$.

\begin{prf} (Proposition 4.13) Note that all but finitely many of the fields $E$ constructed in the preceding lemma are linearly disjoint from $\wti{K}$. Linear disjointness ensures that $[EK:K]=$ $[E:\bq]=2d$. So, from the theory of division algebras, it follows that $EK$ has an embedding in $D$ if and only if $EK$ splits $D$. 

Let $\mf{l}_1,\cdots,\mf{l}_r$ be the non-archimedean primes of $K$ that $D$ is ramified at. Let $l_1,\cdots,l_r$ be the rational primes they lie over. Choose a cyclic extension $E/\bq$ of degree $2d$ inert at the rational primes ${l_1},\cdots,{l_r}$ and linearly disjoint from $K$. Since $\gcd(2d,[\wti{K}:\bq]) =1$, it follows that $EK/K$ is inert at every prime of $K$ lying over $l_i$. In particular, it is inert at every $\mf{l_i}$ and hence, $EK_{\mf{l_i}} = K_{\mf{l_i}^{2d}}$ for every $i$. From the exact sequence \vspace{-0.15cm}$$0\lra \Br(EK|K)\xra{\;\;i\;\;} \bigoplus_w \Br(EK^w|K_w)\xra{\al\mapsto \sum_w \inv_w(\al)} \frac{1}{2d}\bz/\bz\lra 0.$$ we see that $[D\otimes EK_{\mf{l_i}}] = 0$ in $\Br(EK_{\mf{l_i}})$ for each $\mf{l_i}$. Since $EK$ is a CM field, it follows that $[D]$ lies in $\Br(EK|K)$. Hence, $EK$ has an embedding in $D$.\end{prf}

We conclude this article with the following proposition which, when combined with ([No01], Proposition 5.1), allows us to construct abelian varieties strongly of Mumford type such that the reduction at some place is simple.

\begin{Prop} Let $K$ be a totally real number field of degree $N$ with Galois closure $\wti{K}$ over $\bq$ and let $D$ be a quaternion algebra over $K$. Then there exist infinitely many CM fields $L$ quadratic over $K$ such that:\\
- $L$ contains no proper CM subfield\\
- the Galois closure $\wti{L}$ of $L$ over $\bq$ satisfies $\Gal(\wti{L}/\wti{K}) = (\bz/2\bz)^N$\\
- $L$ has an embedding in $D$.\end{Prop}

\begin{prf} Let $\mf{l}_1,\cdots,\mf{l}_r$ be the finite primes of $K$ that $D$ is ramified at and let $l_1,\cdots,l_r$ be the rational primes lying under them. Let $K_1,\cdots, K_t$ be the maximal proper subfields of $K$. We choose an odd rational prime $p$ split completely in $\wti{K}(\sqrt{-1})$ and co-prime to the primes $l_1,\cdots,l_r$. Let $\mfp$ be a prime in $\wti{K}$ lying over $p$. Then the primes $\{\si({\mfp}):\si\in \Gal(\wti{K}/\bq)\}$ are the $[\wti{K}:\bq]$ distinct primes of $\wti{K}$ lying over $p$.

Choose a totally positive $\al\in \mc{O}_K$ such that:\\
- $\bq(\al) = K$.\\
- $X^2+\al$ is irreducible in $(\mc{O}_K/\mf{l}_i)[X]$ for every $i$.\\
- $X^2+\al$ is irreducible in $(\mc{O}_{\wti{K}}/\mfp)[X]$ but splits in $(\mc{O}_{\wti{K}}/\si(\mfp))[X]$ for every $\si\in \Gal(\wti{K}/\bq)-\{\mr{id}\}$.

The existence of the element $\al$ follows from the Chinese remainder theorem. Set $L = K(\sqrt{-\al})$ and let $\wti{L}$ be its Galois closure over $\bq$. Then $\wti{L} = \wti{K}(\{\sqrt{-\si(\al)}:\si\in \Gal(\wti{K}/\bq)\})$. From the second condition, it follows that $L$ has local degree $2$ at each $\mf{l}_i$ and since it is totally imaginary, it follows that $L$ splits $D$. Thus, $L$ has an embedding in $D$. 

For every maximal proper subfield $K_i\sub K$, let $\bar{\mfp}_i$ be a prime in $K_i$ lying under $\mfp$. The third condition implies that not all primes of $L$ lying over $\bar{\mfp}_i$ have the same inertia degree and hence, $L$ is not Galois over $K_i$. Since $K_i$ was an arbitrary maximal subfield of $K$, we have shown that $L$ is not Galois over any proper subfield of $K$. Suppose $L$ has a proper CM subfield $L_0$. Then $L_0\cap K$ is a totally real proper subfield of $K$ of index $2$ in $L_0$. Let $L\cap K \sub K_j$ (notations as above). Then $L = L_0K$ and since $K/\bq$ is Galois, it follows that $L$ is Galois over $L\cap K$ and hence, over $K_j$, a contradiction.

In this paragraph, we slightly abuse the notation by using the same symbols for elements of $\mc{O}_{\wti{K}}$ and their reductions in finite residue fields. The third condition implies that for any distinct elements $\si,\tau\in \Gal(\wti{K}/\bq)$, $X^2+\si(\al)$ is irreducible in $(\mc{O}_{\wti{K}}/\si(\mfp))[X]$ but splits in $(\mc{O}_{\wti{K}}/\tau(\mfp))[X]$. Furthermore, since $p$ splits completely in $\bq(\sqrt{-1})$, we have $-1\in (\mc{O}_{\wti{K}}/\si(\mfp))^{* 2 }$ for every $\si\in \Gal(\wti{K}/\bq)$. Consider an element $\be = \prod\limits_{i=1}^{d}\si_i(\al)$, $1\leq d\leq [\wti{K}:\bq]$. Now $\si_i(\al) \in (\mc{O}_{\wti{K}}/\si_{_1}(\mfp))^{* 2 }$ for every $\si_i\neq \si_{_1}$. On the other hand, $\si_{_1}(\al)\notin (\mc{O}_{\wti{K}}/\si_{_1}(\mfp))^{* 2 }$. Thus, $\be\notin (\mc{O}_{\wti{K}}/\si_{_1}(\mfp))^{* 2 }$ and hence, $\be\notin \wti{K}^{* 2}$. Since the element $\be$ was an arbitrary product of distinct conjugates of $\al$, we have shown that the subgroup of $\wti{K}^*/\wti{K}^{*2}$ generated by the conjugates of $\al$ is of rank $N$. By Kummer theory, it follows that $\Gal(\wti{L}/\wti{K}) = (\bz/2\bz)^{N}$.\end{prf}

Let $N$ be an odd integer, $K$ a totally real number field of degree $N$ and $D$ a quaternion algebra over $K$ ramified at all archimedean places but one such that $\Cor_{K/\bq}(D) = 0$ in $\Br(\bq)$. Let \vspace{-0.15cm}$$\mr{Nm}:D^*\lra \Cor_{K/\bq}(D),\;\;\;\;d\mapsto (d\otimes 1)\otimes \cdots \otimes (d\otimes 1)$$ be the natural norm map, which we may consider as a morphism of algebraic groups over $\bq$. Let $G$ be the image of $D^*$ under this map.

For any quadratic extension $L/K$ constructed as in the preceding proposition, we have $\Gal(\wti{L}/\bq)\cong \{\pm 1\}^{N}\rtimes H$ where $H$ is a transitive subgroup of $S_{_N}$. The map \vspace{-0.15cm}$$h:S\xra{} L_{\br}^*\sub D_{\br}^*\xra{\mr{Nm}}  G_{\br}$$ gives a special point on the Mumford curve which corresponds to an abelian variety $Y$ with CM by the subfield $\wti{L}^H\sub \wti{L}$ fixed by $H$. By [No01], there exists a finite place $v$ of $K$ and a generic point on the curve such that the corresponding abelian variety $X$ has reduction $X_v$ isogenous to $Y_v$. In particular, $X_v$ is simple with CM by the field $\wti{L}^H\sub \wti{L}$.

\bigskip

\textbf{Acknowledgments:} The author is grateful to Ben Howard, David Geraghty and Anand Patel for helpful conversations. We also thank Jeffrey Achter, Anand Patel and Ananth Shankar for helpful suggestions on several previous drafts of the article.

\begin{center}
 \textbf{References} 
 \end{center}
\vspace{-0.2cm}
\footnotesize 

\noindent[Ach12] J. Achter, \textit{Explicit bounds for split reductions of simple abelian varieties}, Journal de Theorie des Nombres de Bordeaux (2012)

\noindent[Bog80] F. Bogomolov, \textit{Sur l'algebricite des representations l-adiques}, C. R. Acad. Sci. Paris S´er. A-B 290 (1980), no. 15.

\noindent[Cha97] N. Chavdarov, \textit{The generic irreducibility of the numerator of the zeta function in a family of curves with large monodromy}, Duke Math. J. 87 (1997), no. 1

\noindent[Chi92] W. Chi, \textit{$l$-adic and $\lambda$-adic representations associated to abelian varieties defined over number fields}, Amer. J. Math. 114 (1992).

\noindent[CCO13] B. Conrad, C. Chai, F. Oort, \textit{Complex multiplication and lifting problems}, AMS series Mathematical Surveys and Monographs, Vol. 195. AMS 2013 Chapter 1.

\noindent[Del74] P. Deligne, \textit{La Conjecture de Weil I}, Publ. Math. IHES 43 (1974)

\noindent[Fal86] G. Faltings, \textit{Finiteness theorems for abelian varieties over number fields}, Arithmetic geometry (Storrs,
Conn., 1984), 1986.

\noindent[FS14] Fukshansky and Sun, \textit{On the geometry of cyclic lattices}, Discrete Comput. Geom., 52(2), 2014
 
\noindent[GRR72] A. Grothendieck, M. Raynaud, and D. S. Rim. \textit{Groupes de monodromie en geometrie algebrique I}, Lecture Notes in Mathematics 288 Springer-Verlag, 1972. Seminaire de Geometrie Algebrique du Bois-Marie 1967-1969

\noindent[Kne65] Kneser, \textit{Galoiskohomologie halbeinfacher algebraischer Gruppen uber p-adischen Korpern}, Teil II.
Math. Zeitschr. 89 (1965).

\noindent[Kow03] E. Kowalski, \textit{Some local-global applications of Kummer theory}, Manuscripta Math., 111(1), 2003.

\noindent[Kow06] \underline{\;\;\;\;\;\;\;} , \textit{Weil numbers generated by other Weil numbers and torsion fields of abelian varieties}, J. London Math. Soc. (2) 74 (2006), no. 2.

\noindent[LP95] M. Larsen and R. Pink, \textit{Abelian varieties, l-adic representations, and l-independence}, Math. Ann. 302 (1995).

\noindent[LP97] \underline{\;\;\;\;\;\;\;} , \textit{A connectedness criterion for l-adic Galois representations}, Israel J. Math. 97 (1997).

\noindent[MP08] V.K. Murty and V. Patankar, \textit{Splitting of abelian varieties}, Int. Math. Res. Not. IMRN 12, 2008

\noindent[Mum69] D. Mumford, \textit{A note on Shimura's paper ``Discontinuous groups and abelian varieties"}, Math. Ann. 181

\noindent[MZ95] B. Moonen and Y. Zarhin, \textit{Hodge classes and Tate classes on simple abelian fourfolds}, Duke Math. J. 77, 3

\noindent[No95] R. Noot, \textit{Abelian varieties- Galois representations and properties of ordinary reduction}, Comp. Math. 97

\noindent[No00] \underline{\;\;\;\;\;\;\;} , \textit{Abelian varieties with l-adic Galois representation of Mumford type}, J. Reine Angew. Math. 519 (2000)

\noindent[No01] \underline{\;\;\;\;\;\;\;} , \textit{On Mumford's families of abelian varieties}, Journal of Pure and Applied Algebra 157 (2001).

\noindent[No09] \underline{\;\;\;\;\;\;\;} , \textit{Classe de conjugaison du Frobenius d’une variete abelienne sur un corps de nombres}, J. Lond. Math. Soc. (2) 79 (2009), no. 1.

\noindent[Orr15] M. Orr, \textit{Lower bounds for ranks of Mumford-Tate groups}, Bull. Soc. Math. France, 143(2), 2015

\noindent[Pin98] R. Pink, \textit{l-adic algebraic monodromy groups, co-characters and the Mumford-Tate conjecture}, J. Reine Angew. Math. 495 (1998)

\noindent[Ser79] J.P. Serre, \textit{Groupes algebriques associes aux modules de Hodge-Tate}, Journees de Geometrie Algebrique de Rennes. (Rennes, 1978), Vol. III, 1979

\noindent[Ser00] \underline{\;\;\;\;\;\;\;} , \textit{OEuvres. Collected papers. IV}, Springer-Verlag, Berlin, 2000.

\noindent[Ser03] \underline{\;\;\;\;\;\;\;} , \textit{On a theorem of Jordan}, Bull. Amer. Math. Soc. 40 (2003).

\noindent[ST68] J.P. Serre and J.Tate, \textit{Good reduction of abelian varieties}, Ann. of Math. (2), 88, 1968

\noindent[SZ12] M. Sheng and K. Zuo, \textit{On the Newton polygons of abelian varieties of Mumford type}, Preprint

\noindent[Tat66] J. Tate, \textit{Endomorphisms of abelian varieties over finite fields}, Invent. Math, Vol 2, 1966

\noindent[Wat69] W. Waterhouse, \textit{Abelian varieties over finite fields}, Ann. Sci. Ecole Norm. Sup. (4) 2 (1969).

\noindent[Win02] J.P. Wintenberger, \textit{Demonstration d’une conjecture de Lang dans des cas particuliers}, J. Reine Angew. Math. 553 (2002).

\noindent[Zyw14] D. Zywina, \textit{The splitting of reductions of an abelian variety}, IMRN, vol. 2014, No. 18.

\bigskip

\noindent E-mail address: steve.thakur@bc.edu, stevethakur01@gmail.com\\
Mathematics Department\\
Boston College.

\end{document}